\documentclass[10pt,draft
]{article}

\usepackage[latin1]{inputenc}
\usepackage{amsfonts}
\usepackage{amsmath}
\usepackage{amssymb}
\usepackage{amscd}
\usepackage{amsthm}
\usepackage{indentfirst}
\usepackage[hmargin=3.6cm
,vmargin=4.1cm
]{geometry}

\usepackage{epsfig}

\newtheorem{thm}{Theorem}[section]
\newtheorem{lemma}[thm]{Lemma}
\newtheorem{prop}[thm]{Proposition}
\newtheorem{coroll}[thm]{Corollary}

\theoremstyle{definition}

\newtheorem{defin}[thm]{Definition}
\newtheorem{rem}[thm]{Remark}
\newtheorem{exam}[thm]{Example}

\newtheorem{notation}[thm]{Notation}
\newtheorem*{acknow}{Acknowledgements}
\newtheorem*{prf}{Proof}

\newcommand{\R}{{\mathbb{R}}}

\newcommand{\T}{{\mathbb{T}}}
\newcommand{\Z}{{\mathbb{Z}}}
\newcommand{\N}{{\mathbb{N}}}
\newcommand{\C}{{\mathbb{C}}}

\newcommand{\cA}{{\mathcal{A}}}

\newcommand{\cG}{{\mathcal{G}}}
\newcommand{\cH}{{\mathcal{H}}}

\newcommand{\cM}{{\mathcal{M}}}

\newcommand{\cO}{{\mathcal{O}}}
\newcommand{\cP}{{\mathcal{P}}}

\newcommand{\cS}{{\mathcal{S}}}

\newcommand{\cU}{{\mathcal{U}}}
\newcommand{\cV}{{\mathcal{V}}}

\newcommand{\fc}{{:\ }}
\newcommand{\ve}{\varepsilon}

\newcommand{\ol}{\overline}

\newcommand{\tb}{\textbf}

\DeclareMathOperator{\Hom}{Hom}

\DeclareMathOperator{\Crit}{Crit}

\DeclareMathOperator{\im}{im}
\DeclareMathOperator{\id}{id}

\DeclareMathOperator{\osc}{osc}

\DeclareMathOperator{\End}{End}
\DeclareMathOperator{\Spec}{Spec}

\title{Partial quasi-morphisms and quasi-states on cotangent bundles, and symplectic homogenization}
\author{Alexandra Monzner,\footnote{Fakult\"at f\"ur Mathematik, TU Dortmund, Dortmund, Germany,\newline \texttt{alexandra.monzner@mathematik.tu-dortmund.de}}\; Nicolas Vichery,\footnote{CMLS \'Ecole Polytechnique, Palaiseau, France, \texttt{nicolas.vichery@math.polytechnique.fr}}\; and Frol Zapolsky\footnote{IHES, Bures-sur-Yvette, France, \texttt{zapolsky@ihes.fr}}}
\date{}

\begin{document}

\renewcommand{\labelenumi}{(\roman{enumi})}

\maketitle

\begin{abstract}
For a closed connected manifold $N$, we construct a family of functions on the Hamiltonian group $\cG$ of the cotangent bundle $T^*N$, and a family of functions on the space of smooth functions with compact support on $T^*N$. These satisfy properties analogous to those of partial quasi-morphisms and quasi-states of Entov and Polterovich. The families are parametrized by the first real cohomology of $N$. In the case $N=\T^n$ the family of functions on $\cG$ coincides with Viterbo's symplectic homogenization operator. These functions have applications to the algebraic and geometric structure of $\cG$, to Aubry-Mather theory, to restrictions on Poisson brackets, and to symplectic rigidity.
\end{abstract}

\tableofcontents

\section{Introduction and results}\label{section_intro_and_results}

\subsection{Overview}\label{section_overview_of_results}

Fix a closed connected manifold $N$ of dimension $n$. The cotangent bundle $T^*N$ has a natural symplectic structure. We let $\cG$ be the Hamiltonian group with compact support of $T^*N$. We construct two families of functions, $\mu_a \fc \cG \to \R$, and $\zeta_a \fc C^\infty_c(T^*N) \to \R$, where $a \in H^1(N;\R)$. These functions possess properties analogous to those of partial quasi-morphisms and partial quasi-states of Entov and Polterovich \cite{EP_qs_sympl}.\footnote{Strictly speaking, the term ``partial quasi-morphism'' is not attested in the existing literature, as far as we know. But it is known and used in the community, so there is no harm in utilizing this term here, which will also serve to its dissemination.} The precise properties are listed in theorems \ref{thm_main_result}, \ref{thm_properties_zeta} below.

In the case $N=\T^n$, the family $\mu_a$ is equivalent to Viterbo's symplectic homogenization. The symplectic homogenization is an operator (see \cite{Viterbo_homogenization})
$$C^\infty_c([0,1]\times T^*\T^n) \to C_c(\R^n)\,,\quad H \mapsto\ol H\,.$$
Identify $H^1(\T^n;\R)=\R^n$. Then we have
\begin{thm}\label{thm_our_constr_equals_homogenization}Let $N = \T^n$. Then $\ol H(p)$ equals the value of $\mu_p$ on the time-$1$ map of $H$, for any $H \in C_c^\infty([0,1]\times T^*\T^n)$ and any $p \in \R^n$.
\end{thm}

The properties of $\mu_a,\zeta_a$ lead to various applications. Briefly, these include lower bounds on the fragmentation norm on $\cG$ relative to displaceable subsets, symplectic invariance of Mather's alpha function, Hofer and spectral geometry of $\cG$, restrictions on the Poisson brackets and symplectic rigidity of subsets of $T^*N$. Most of these applications have appeared in the literature in some form or another; we indicate the connection to the existing results where appropriate. We would like to point out, however, that here we present a unified approach which provides transparent and elementary proofs of all of the above results, together with new ones.

The main technical tool in the construction of $\mu_a$ and $\zeta_a$ is the spectral invariants in Lagrangian Floer homology, which themselves are functions on $\cG$. The $\mu_a$ are obtained from the spectral invariants via homogenization, and $\zeta_a$ are obtained from $\mu_a$ by pulling them back via the exponential map.

Spectral invariants have been used for some time now to prove interesting and deep results in symplectic topology; to list but a few references: \cite{Viterbo_gfqi}, \cite{Oh_action_I}, \cite{EP_Calabi_qm}. Our contribution to their theory in this paper is twofold. Firstly, we prove a sharp triangle inequality for them, which implies, in particular, that the invariants descend to the Hamiltonian group. Secondly, we prove an inequality relating spectral invariants coming from Lagrangian and Hamiltonian Floer homology. This allows us to prove a vanishing property for the $\mu_a$ and $\zeta_a$.

Floer-homological spectral invariants have become standard objects in symplectic topology, and in our opinion the fact that symplectic homogenization is expressible with their help, makes the latter fit nicely into the general theory.

The rest of the paper is organized as follows. The remainder of this section is devoted to precise formulations of the properties of the $\mu_a$ and $\zeta_a$ and their applications. In section \ref{section_spectral_invariants} we present the construction and properties of Lagrangian and Hamiltonian spectral invariants on $\cG$. Section \ref{section_proofs} contains proofs of the results formulated in subsections \ref{section_main_result} and \ref{section_Applications}. The reader interested in proofs of the applications can go directly to section \ref{section_proofs}, after reviewing subsection \ref{section_summary_sp_invts}. With rare exceptions, the proofs presented in section \ref{section_proofs} rely only on the properties of the spectral invariants appearing there.

\subsubsection{Preliminaries and notations}\label{section_conventions_notations}

The symplectic form on $T^*N$ is $\omega = d\lambda = dp\wedge dq$, where $\lambda = p\,dq$ is the Liouville form. The zero section of a cotangent bundle $T^*Q$ is denoted by $Q$, unless a confusion may arise, in which case we employ the more explicit notation $\cO_Q$.

We implicitly fix an auxiliary Riemannian metric on $N$ and other closed manifolds appearing below, and the lengths of cotangent vectors are measured relative to this metric.

A time-dependent Hamiltonian, that is a smooth function $H \fc [0,1]\times T^*N \to \R$, is either denoted by $H$ or by explicitly pointing out the time-dependence, $H_t$. This symbol also means the function $H(t,\cdot) \in C^\infty(T^*N)$. The time-$t$ map of the flow of $H$ is denoted by $\phi_H^t$ and the time-$1$ map by $\phi_H$. The collection of time-$1$ maps of all the Hamiltonians with compact support is the Hamiltonian group $\cG$ of $T^*N$.

For an open subset $U \subset T^*N$ we let $\cG_U \subset \cG$ be the subgroup generated by Hamiltonians with compact support in $U$. We let $T^*_{<r}N = \{(q,p)\,|\,\|p\|<r\}$, for $r > 0$.

An interesting subgroup of $\cG$ consists of all the Hamiltonian diffeomorphisms fixing the zero section $N$ as a set. It is denoted by $\cG_0$.
\begin{prop}\label{prop_action_homomorphism}There is a natural action homomorphism $\cA \fc \cG_0 \to \R$.
\end{prop}
\noindent See subsection \ref{section_arbitrary_Hamiltonian} for a precise formulation and a proof. For now let us just note that if $H$ is a time-dependent Hamiltonian which equals $c \in \R$ when restricted to the zero section, then $\cA(\phi_H)=c$.

A subset $S \subset T^*N$ is called displaceable if there is $\phi \in \cG$ with $\ol S \cap \phi(S)=\varnothing$. We say that $S$ is dominated by an open subset $U$ if there is $\phi \in \cG$ such that $S \subset \phi(U)$. In subsection \ref{section_Ham_sp_invts} we introduce the spectral norm $\Gamma \fc \cG \to \R$. The spectral displacement energy of a displaceable subset $S$ is by definition $e(S)=\inf\{\Gamma(\psi)\,|\,\psi(S)\cap\ol S = \varnothing\}$. The spectral displacement energy of a family $\cS = \{S_i\}_i$ of subsets is $e(\cS)=\sup_ie(S_i)$.

We also introduce the fragmentation norm. This is defined as follows. If $\cV$ is an open covering of $T^*N$, then Banyaga's fragmentation lemma \cite{Banyaga_Ham} states that any $\phi \in \cG$ can be represented as a finite product $\phi = \prod_i \phi_i$ where every $\phi_i$ belongs to $\cG_{U_i}$ for some $U_i \in \cV$. The fragmentation norm of $\phi$ relative to the covering $\cV$ is the minimal number of such factors. We will need the following version of the fragmentation norm. Let $\cU$ be an arbitrary family of open subsets and consider the open covering $\cV$ consisting of all open subsets $V$ for which there is $\psi \in \cG$ with $\psi(V) \in \cU$. We let $\|\phi\|_\cU$ be the fragmentation norm of $\phi$ relative to the covering $\cV$.

We define $\phi_H^t := \phi_H^{t - k}\phi_H^k$ for $t \in [k,k+1]$, where $k \in \Z$; here $\phi_H^k := (\phi_H)^k$. Whenever $H$ is defined for all $t\in\R$ and is $1$-periodic in $t$, the time-$t$ flow of $H$ equals $\phi^t_H$.

\subsection{Properties of $\mu_a$ and $\zeta_a$}\label{section_main_result}

The following theorem lists the properties of $\mu_a$. Recall the notion of the spectral displacement energy $e$ of a family of subsets introduced above.

\begin{thm}\label{thm_main_result}Let $N$ be a closed connected manifold. For every $a \in H^1(N;\R)$ there is a function $\mu_a \fc \cG \to \R$ with the following properties:
\begin{enumerate}
\item $\mu_a(\phi^k)=k\mu_a(\phi)$ for $k \geq 0$ an integer;
\item $\mu_a$ is conjugation-invariant;
\item if $\phi,\psi\in\cG$ are generated by the Hamiltonians $H,K$, then
$$\int_0^1\min(H_t-K_t)\,dt \leq \mu_a(\phi)-\mu_a(\psi) \leq \int_0^1\max(H_t-K_t)\,dt\,;$$
in particular $\mu_a$ is Lipschitz with respect to the Hofer metric;
\item the restriction of $\mu_a$ to $\cG_U$ vanishes for any displaceable $U$;
\item for any collection $\cU$ of open subsets with $e(\cU) < \infty$ we have
$$|\mu_a(\phi\psi)-\mu_a(\psi)| \leq e(\cU)\|\phi\|_\cU\,;$$
\item the restriction of $\mu_0$ to $\cG_0$ coincides with the action homomorphism $\cA$;
\item if $\phi \in \cG$ is generated by a Hamiltonian whose restriction to the graph of a closed $1$-form in the class $a$ is $\geq c$ (respectively, $\leq c$, $=c$), where $c$ is some number, then
$$\mu_a(\phi) \geq c\;(\text{respectively}\leq c,=c)\,;$$
\item for commuting $\phi,\psi$ we have $\mu_a(\phi\psi) \leq \mu_a(\phi)+\mu_a(\psi)$;
\item for fixed $\phi \in \cG$ the function $H^1(N;\R)\to\R$, $a \mapsto \mu_a(\phi)$, is Lipschitz, the Lipschitz constant being given by a semi-norm.
\end{enumerate}
\end{thm}
\noindent We may call these $\mu_a$ partial quasi-morphisms in the sense of Entov-Polterovich.

\begin{rem}Combining this theorem with theorem \ref{thm_our_constr_equals_homogenization} we see that now we have a definition of the symplectic homogenization (as an operator $\cG \to C_c(H^1(N;\R)$) for any base. In fact, the proof of convergence in \cite{Viterbo_homogenization} relies on more assumptions than that of the existence of $\mu_a$.\footnote{It relies, as far as we understand, on the existence of certain capacity bounds, see also subsection \ref{section_generalizations} below.} This means that theorem \ref{thm_main_result} gives an alternative \emph{definition} of symplectic homogenization. The properties of $\mu_a$, listed in this theorem, give properties of symplectic homogenization; in particular, the Lipschitz property of $\ol H$ mentioned \emph{ibid.} is a consequence of point (ix) of our theorem.
\end{rem}

In the original work \cite{Viterbo_homogenization} the author constructs the symplectic homogenization as a limit in certain variational metric of a sequence of flows when one passes to coverings of arbitrary large degrees. More precisely, consider the conformal symplectic covering $r_k \fc T^*\T^n \to T^*\T^n$, $r_k(q,p)=(kq,p)$. Hamiltonian flows can be pulled back via this covering, namely if $H$ is a time-dependent Hamiltonian generating $\phi$, put $H_k(t,q,p)=H(kt,kq,p)$ and let $\phi_k$ be the time-$1$ map of $H_k$. Then the symplectic homogenization of $\phi$ is a continuous Hamiltonian which only depends on $p$, which generates in a certain precise sense the limit in the aforementioned metric of the sequence $\phi_k$. At least philosophically, it follows from the equivalence of the symplectic homogenization and our functionals $\mu_p$, $p\in \R^n$, that the latter are invariant under this passage to coverings. More precisely, we have the following claim.
\begin{prop}\label{prop_mu_p_invt_under_coverings}For any $k$ we have $\mu_p(\phi_k) = \mu_p(\phi)$, $\phi \in \cG$, $\phi_k$ being defined as above. \qed
\end{prop}
\noindent This can be extracted as a byproduct of the proof of theorem \ref{thm_our_constr_equals_homogenization}.
\begin{rem}As pointed out by L.\ Polterovich, this result makes our construction fit nicely in the philosophy of homogenization, which in particular manifests itself in such objects as the Gromov-Federer stable norm, and on the other hand shows that the classical notion of homogenization, in this sense, is a particular case of a Floer-homological construction applicable to general (not necessarily convex) Hamiltonian dynamical systems.
\end{rem}

For applications it is important to extend the definition of $\mu_a$ to more general diffeomorphisms. In subsection \ref{section_geom_bdd_Hamiltonians} we show how to define $\mu_a(\phi_H)$ in case $H$ is a Hamiltonian with complete flow. For time-dependent Hamiltonians $H_t,H_t'$ on symplectic manifolds $Z,Z'$, respectively, we define the direct sum $H \oplus H'$ to be the time-dependent Hamiltonian on $Z\times Z'$ given by $(H\oplus H')(t,z,z')=H(t,z)+H'(t,z')$. An easy observation is that whenever $H,H'$ have complete flows, so does their sum $H \oplus H'$. With this observation at hand we formulate
\begin{prop}\label{prop_product_for_qm}Assume that $N=N_1\times N_2$, and that $\mu_{a_i}^{(i)}$, $i=1,2$ are the corresponding functionals given by theorem \ref{thm_main_result}, extended to the set of complete flows, where $a_i \in H^1(N_i;\R)$. Let $H^{(i)}$ be a time-dependent Hamiltonian on $T^*N_i$ for $i=1,2$, and assume both have complete flows. Then for $a = (a_1,a_2) \in H^1(N_1;\R)\times H^1(N_2;\R)\subset H^1(N;\R)$ we have
$$\mu_a(\phi_{H^{(1)}\oplus H^{(2)}}) = \mu_{a_1}^{(1)}(\phi_{H^{(1)}})+\mu_{a_2}^{(2)}(\phi_{H^{(2)}})\,.$$
\end{prop}

We let $\zeta_a \fc C^\infty_c(T^*N)\to \R$ be defined as $\zeta_a(H) = \mu_a(\phi_H)$. The following theorem lists the properties of $\zeta_a$.
\begin{thm}\label{thm_properties_zeta}
\begin{enumerate}
\item $\zeta_a(\lambda F)=\lambda \zeta_a(F)$ for $\lambda \geq 0$ a real number;
\item $\zeta_a$ is invariant under the natural action of $\cG$ on $C^\infty_c(T^*N)$;
\item $\min(F-G) \leq \zeta_a(F)-\zeta_a(G) \leq \max (F-G)$, in particular $|\zeta_a(F)-\zeta_a(G)| \leq \|F-G\|_{C^0}$;
\item $\zeta_a(F)=0$ for $F$ with displaceable support;
\item for displaceable $U$, any $F\in C^\infty_c(T^*N)$ and any $G$ with support dominated by $U$, we have
$$|\zeta_a(F+G) - \zeta_a(F)-\zeta_a(G)| \leq \sqrt{2e(U)\|\{F,G\}\|_{C^0}}\,;$$
in particular, if $F,G$ commute and the support of $G$ is displaceable then $\zeta_a(F+G)=\zeta_a(F)+\zeta_a(G)=\zeta_a(F)$;
\item if $F \geq c$ (respectively, $\leq c$) when restricted to the graph of a closed $1$-form in the class $a$, then $\zeta_a(F)\geq c$ (respectively, $\leq c$);
\item if $\{F,G\}=0$ then $\zeta_a(F+G)\leq\zeta_a(F)+\zeta_a(G)$.
\end{enumerate}
\end{thm}
Similarly to $\mu_a$, $\zeta_a$ can be defined on autonomous Hamiltonians with complete flow. For these we have the following product formula, which follows from the one formulated in proposition \ref{prop_product_for_qm}:
\begin{prop}Assume that $N=N_1\times N_2$ and that $\zeta_{a_i}^{(i)}$, $i=1,2$, are the corresponding functionals; then if $F_i \in C^\infty(T^*N_i)$ have complete flows, then for $a=(a_1,a_2) \in H^1(N_1;\R)\times H^1(N_2;\R) \subset H^1(N;\R)$
$$\zeta_a(F_1 \oplus F_2)=\zeta_{a_1}^{(1)}(F_1)+\zeta_{a_2}^{(2)}(F_2)\,. \qed$$
\end{prop}

\subsection{Applications}\label{section_Applications}

\subsubsection{Fragmentation norm}\label{section_frag_comm_norm}

\begin{coroll}\label{coroll_frag_comm_norm}
The fragmentation norm relative to a family of open subsets $\cU$ satisfies:
$$\|\phi\|_{\cU} \geq \sup_{a\in H^1(N;\R)}\frac{|\mu_a(\phi)|}{e(\cU)}\,.$$
In particular, if $\phi$ is generated by a Hamiltonian whose restriction to $L$ is at least $c$ in absolute value, where $L$ is a Lagrangian submanifold Hamiltonian isotopic to the zero section, and $c$ is a number, then
$$\|\phi\|_\cU \geq \frac c {e(\cU)}\,.$$
\end{coroll}
\begin{prf}Point (v) of theorem \ref{thm_main_result} implies (with $\psi = \id$) that
$$|\mu_a(\phi)| \leq e(\cU)\|\phi\|_\cU\,.$$
For the second claim it suffices to note that for such $\phi$ we have $|\mu_0(\phi)|\geq c$. \qed
\end{prf}

\noindent Similar results are proved in \cite{EP_Calabi_qm}, \cite{Lanzat_qms_convex_mfds}. The difference is in the class of manifolds under consideration (closed manifolds in the first reference and certain types of open convex manifolds in the second, including the unit disk cotangent bundle of a torus) and, in the case of the first reference, that there the Hamiltonian diffeomorphism $\phi$ is itself required to have displaceable support. This has to do with the fact that the Calabi quasi-morphism used there coincides with the Calabi invariant on displaceable subsets while our $\mu_a$ (and Lanzat's functionals) vanish on displaceable subsets.

\subsubsection{Connection with Mather's alpha function}\label{section_connection_alpha_fcn}

Aubry-Mather theory, among other things, associates a function on $H^1(N;\R)$ to a Tonelli Hamiltonian $H$, the so-called alpha function $\alpha_H \fc H^1(N;\R) \to \R$. A Hamiltonian $H \fc [0,1]\times T^*N \to\R$ is called Tonelli if it is fiberwise strictly convex and superlinear, and has complete flow. We refer the reader to \cite{Mather_action_minimizing} for a more detailed exposition. The functions $\mu_a$ appearing in theorem \ref{thm_main_result} can be correctly defined on Hamiltonians having complete flow. This is done in subsection \ref{section_geom_bdd_Hamiltonians}. We have
\begin{thm}\label{thm_homogen_equals_alpha}Let $H$ be a time-periodic Tonelli Hamiltonian. Then for $a \in H^1(N;\R)$
$$\alpha_H(a)=\mu_a(\phi_H)\,.$$
\end{thm}
\noindent One way to interpret this result is that now we have a way of defining the alpha function for an arbitrary Hamiltonian $H$ with complete flow: $\alpha_H(a):=\mu_a(\phi_H)$. Theorem \ref{thm_homogen_equals_alpha} first appeared in \cite{Viterbo_homogenization} in the case $N =\T^n$.

An immediate consequence of this theorem is formulated in the following
\begin{coroll}\label{coroll_props_of_alpha_fcn}Let $H$ be a time-periodic Tonelli Hamiltonian. Then if $\phi \in \cG$ is such that $H\circ \phi$ is still Tonelli,\footnote{Since $\phi$ has compact support, $H\circ\phi$ automatically has complete flow.} then
$$\alpha_{H\circ\phi} = \alpha_H\,.$$
\end{coroll}
\begin{prf}The extended functionals $\mu_a$ are still invariant under conjugation by elements of $\cG$. Since $H\circ \phi$ generates the diffeomorphism $\phi^{-1}\phi_H\phi$, the desired conclusion follows from this conjugation invariance. \qed
\end{prf}
\noindent The reader can find the proof of the symplectic invariance of the alpha function in \cite{Bernard_sympl_aspects} and the references therein, in the case of Tonelli Hamiltonians. It is also implicit in \cite{Paternain_P_S_boundary_rigidity}, in case $H$ is autonomous. The advantage of our approach is that this invariance follows formally from the conjugation invariance of $\mu_a$, and it is applicable to any Hamiltonian with complete flow.

\subsubsection{Hofer geometry and spectral norm on $\cG$}\label{section_Hofer_spectral_norms}

For $\phi \in \cG$ put
$$\rho(\phi) = \inf_H\int_0^1\osc H_t\,dt\,,$$
where $\osc = \max - \min$ and the infimum is over all the compactly supported Hamiltonians whose time-$1$ map is $\phi$. Also put
$$\rho(\phi,\psi) = \rho(\phi\psi^{-1})\,.$$
It is a highly nontrivial fact that $\rho$ is a metric on $\cG$, called the Hofer metric. It is biinvariant.\footnote{The reader is referred to \cite{Polterovich_geom_grp_sympl_diffeo} for preliminaries on Hofer geometry.}

There is another norm on $\cG$, various variants of which were introduced by Viterbo, Schwarz, Oh, and in the present context, by Frauenfelder and Schlenk \cite{Frauenfelder_Schlenk_Ham_convex}. Namely, there are two spectral invariants $c_\pm \fc \cG \to \R$ and the spectral norm is defined to be
$$\Gamma(\phi) = c_+(\phi) - c_-(\phi)\,.$$
See subsection \ref{section_Ham_sp_invts} for more details. Since this norm is conjugation-invariant \cite{Frauenfelder_Schlenk_Ham_convex}, it gives rise to another biinvariant metric on $\cG$, which we call the spectral metric, via
$$\Gamma(\phi,\psi)=\Gamma(\phi\psi^{-1})\,.$$
It is known \cite{Frauenfelder_Schlenk_Ham_convex} that
$$\Gamma(\phi,\psi)\leq\rho(\phi,\psi)\,.$$

For the next theorem, note that oscillation is a norm on the space $C^\infty_c(0,1)$. We denote by $(C^\infty_c(0,1),\osc)$ the corresponding metric space.

\begin{thm}\label{thm_Hofer_embeddings}(i) If $N$ admits a non-singular closed $1$-form, then there are isometric embeddings $(C^\infty_c(0,1),\osc)$ into $(\cG,\rho)$; precisely, there are maps $\iota \fc C^\infty_c(0,1) \to \cG$, such that
$$\rho(\iota(f),\iota(g)) = \osc(f-g)\,;$$
(ii) otherwise there is an isometric embedding of $\R$ into $\cG$;
(iii) the same holds if we replace the Hofer metric with the spectral metric.
\end{thm}
\noindent Contrast this with \cite{Py_plats_pour_Hofer}, where the author constructs, using the energy-capacity inequality, quasi-isometric embeddings of $\R^k$, $k \geq 1$, into the Hamiltonian group of a symplectic manifold admitting a $\pi_1$-injective Lagrangian embedding of a Riemannian manifold of non-positive sectional curvature.

We define the asymptotic Hofer norm
$$\rho_\infty(\phi)=\lim_{k \to \infty}\frac{\rho(\phi^k)}{k}\,.$$
As with the Hofer norm, we can introduce the asymptotic version
$$\Gamma_\infty(\phi) = \lim_{k\to\infty}\frac{\Gamma(\phi^k)}{k}\,.$$
We then have
\begin{prop}\label{prop_asym_Hofer_spectral}Let $\phi \in \cG$. Then
$$\osc_{a\in H^1(N;\R)}\mu_a(\phi) \leq \Gamma(\phi)\leq\rho(\phi)\,;$$
homogenizing, we obtain
$$\osc_{a\in H^1(N;\R)}\mu_a(\phi) \leq \Gamma_\infty(\phi)\leq\rho_\infty(\phi)\,.$$
\end{prop}
\noindent Related results can be found in \cite{Polterovich_Siburg_asymp_geom_area}, \cite{Siburg_book}, \cite{Sorrentino_Viterbo_act_min_pties_dist_Ham}, \cite{Monzner_Zapolsky_comparison}.

There is also a connection between Aubry-Mather theory and Hofer geometry, as studied in \cite{Siburg_act_min_meas_geom_Ham}. We let $\cH$ be the space of Hamiltonian functions on the closed unit disk cotangent bundle $B \subset T^*N$ which vanish at the boundary and which admit smooth extensions to the whole cotangent bundle which only depend on $\|p\|$ and $t$ outside the unit ball bundle. There is the associated notion of Hofer norm:
$$\rho_\cH(\phi)=\inf_{H\in\cH}\int_0^1\osc H_t\,dt\,,$$
where $\phi \fc B \to B$ is the time-$1$ map of a Hamiltonian in $\cH$ and $H$ runs over all Hamiltonians in $\cH$ generating $\phi$. We have
\begin{thm}\label{thm_alpha_fcn_Hofer_metric}Let $\widetilde H$ be a Tonelli Hamiltonian which vanishes for $\|p\|=1$ and which only depends on $\|p\|$ for $\|p\|\geq1$. Let $H=\widetilde H|_B \in \cH$. Then
$$\rho_\cH(\phi_H) \geq -\min_{H^1(N;\R)}\alpha_{\widetilde H}\,.$$
\end{thm}
\noindent This was proved in \cite{Siburg_act_min_meas_geom_Ham} for $N=\T^n$ and in \cite{Iturriaga_Morgado_minimax} for a class of Hamiltonians on the cotangent bundles over a general base, using different methods. Note that the minimum in the right-hand side only depends on $H$. Of course, since we have a definition of the alpha function for any Hamiltonian with complete flow, and the Hofer norm is defined for any compactly supported Hamiltonian, proposition \ref{prop_asym_Hofer_spectral} provides a more natural formulation of the relation between the Hofer norm and the alpha function, so we only include this result for completeness's sake and to illustrate the power of the methods developed here.

\subsubsection{Poisson brackets and symplectic rigidity}\label{section_qi_symp_rigidity}

We abbreviate $\zeta=\zeta_0$. Property (v) of $\zeta$ implies the following restrictions on Poisson brackets.
\begin{thm}There are constants\footnote{It is true that $C\geq \frac 9 8$ and $C'\geq 1/2$.} $C,C' > 0$ such that the following holds. If $\{f_i\}_{i=1}^K$ are smooth functions such that the support of each one of them is dominated by an element in a fixed collection $\cU$ of displaceable subsets, and which satisfy $\left.\sum_if_i\right|_{N} \geq 1$, then
$$\max_{i<j}\|\{f_i,f_j\}\|_{C^0} \geq \frac C {e(\cU)K^3}\,.$$
Moreover, if there is a number $k$ such that the number of supports of the $f_j$ intersecting at any point of $T^*N$ is at most $k$, then
$$\max_{i<j}\|\{f_i,f_j\}\|_{C^0} \geq \frac {C'} {e(\cU)kK^2}\,.\qed$$
\end{thm}
\noindent The proof is a \emph{verbatim} repetition of the one in \cite{EPZ_qm_Poisson_br} and will be omitted.

We now turn to non-displaceability. We refer the reader to \cite{EP_qs_sympl}, \cite{EP_rigid_subsets} for a treatment of the rigidity of subsets in closed symplectic manifolds.

Following \cite{EP_rigid_subsets}, we make
\begin{defin}Call a compact subset $X\subset T^*N$ $\zeta$-superheavy, or superheavy for brevity, if for any function $f \in C^\infty_c(T^*N)$ with $f|_X = c \in \R$ we have $\zeta(f) = c$.
\end{defin}
\begin{rem}Since $\zeta$ is invariant under the action of $\cG$, so is the collection of superheavy subsets. Also we would like to point out that in this paper only superheavy (not heavy) subsets appear, since it is easy to construct examples of superheavy subsets but we could not find a heavy subset which is not superheavy.
\end{rem}

\begin{exam}The zero section is superheavy by property (vi) in theorem \ref{thm_properties_zeta}.
\end{exam}

\begin{lemma}\label{lemma_visible_iff_superheavy}A subset $X$ is superheavy if and only if for any $f$ we have $\zeta(f) \leq \max_X f$.
\end{lemma}
\noindent Lemma \ref{lemma_visible_iff_superheavy} is proved in subsection \ref{section_proofs_symp_rigidity}. In fact, the original definition of a superheavy subset used this weaker characterization, which is more easily checked.

Superheavy subsets are rigid in the sense that any two must intersect:
\begin{prop}\label{prop_visible_intersect}Let $X,X'$ be two superheavy subsets; then $X \cap X' \neq \varnothing$.
\end{prop}
\noindent This implies in particular that superheavy subsets are non-displaceable. The proof is short and instructive, thus we include it here.
\begin{prf}Assume the contrary and choose $f,f' \in C^\infty_c(T^*N)$ such that $f|_X = f'|_{X'}=-1$ and $\|f\|_{C^0}=\|f'\|_{C^0} =1$, and such that the supports of $f,f'$ are disjoint. In particular this means that they Poisson commute. Then we have, by property (vii) of $\zeta$
$$\zeta(f+f') \leq \zeta(f)+\zeta(f') = -2\,,$$
which is a contradiction to $|\zeta(f+f')| \leq \|f+f'\|_{C^0}=1$. \qed
\end{prf}
\noindent Using the same argument, one can show that if $X$ is superheavy and has a finite number of connected components, then only one of these connected components is superheavy.

The following proposition, proved in subsection \ref{section_proofs_symp_rigidity} allows us to construct many examples of superheavy subsets.
\begin{prop}\label{prop_examples_of_visible_sets}Let $X$ be a compact subset such that $T^*N-X = U_\infty \cup \bigcup_iU_i$ is a finite disjoint union with $U_\infty$ being the unbounded connected component (the union of the unbounded connected components in case $\dim N = 1$). Assume that $U_\infty$ is disjoint from the zero section and that each one of $U_i$ is displaceable. Then $X$ is superheavy.
\end{prop}

\begin{exam}The codimension $1$ skeleton of a triangulation (or, more generally, a polygonal subdivision) of the closed unit disk cotangent bundle in $T^*N$, considered as a manifold with boundary, satisfies the assumptions of the proposition and thus is superheavy.
\end{exam}

Finally, in order to obtain yet more examples, we formulate the following result, also proved in subsection \ref{section_proofs_symp_rigidity}.
\begin{thm}\label{thm_product_visible}Let $X_i \subset T^*N_i$, $i=1,2$, be superheavy subsets; then the product $X_1\times X_2 \subset T^*N_1 \times T^*N_2 = T^*(N_1\times N_2)$ is superheavy.
\end{thm}

\noindent This implies the following
\begin{coroll}Let $X_i,X_i' \subset T^*N_i$, $i=1,\dots,k$ be subsets as in proposition \ref{prop_examples_of_visible_sets}. Then $\prod_iX_i \cap \phi\big(\prod_iX_i'\big) \neq \varnothing$ for any Hamiltonian diffeomorphism $\phi$ on $T^*\prod_iN_i$. In particular, $\prod_iX_i$ is non-displaceable. \qed
\end{coroll}

\subsection{Connection with existing constructions and generalizations}\label{section_generalizations}

Here we indicate connections to analogous constructions. This is an expository subsection, therefore no proofs are given.

Lanzat \cite{Lanzat_thesis}, \cite{Lanzat_qms_convex_mfds} produces examples of open symplectic manifolds whose Hamiltonian group with compact support (or its universal cover) admits genuine (not partial) quasi-mor\-phisms, and whose space of smooth functions with compact support, as a consequence, admits a symplectic quasi-state. He also shows how to construct a partial quasi-morphism and a partial symplectic quasi-state on a general (strongly semipositive) convex manifold. In particular his construction applies to cotangent bundles. The spectral invariant $c_+\fc \cG \to \R$, introduced in subsection \ref{section_Ham_sp_invts}, satisfies the triangle inequality and so can be homogenized to yield a functional $\nu \fc \cG \to \R$, whose pullback to $C_c^\infty(T^*N)$ is denoted by $\eta$. These $\nu$, $\eta$, in fact, coincide with Lanzat's functionals for cotangent bundles, and enjoy properties analogous to those of $\mu_a$, $\zeta_a$. Moreover, owing to the comparison of Lagrangian and Hamiltonian spectral invariants (subsection \ref{section_comparison_Ham_Lagr_sp_invts}), we can conclude that $\mu_a \leq \nu$ and $\zeta_a \leq \eta$ for any $a \in H^1(N;\R)$. In particular, any $\zeta_a$-superheavy set is $\eta$-heavy and so sets described in proposition \ref{prop_examples_of_visible_sets}. This means that $\eta$ can be used to prove non-displaceability of such subsets. However, since $\eta(f)=0$ for nonpositive functions $f$, the collection of $\eta$-superheavy sets is empty, and therefore the applications to symplectic rigidity end there. It is instructive to note here that in contrast, the collection of $\zeta_a$-superheavy subsets is not empty, and this allows for more flexible rigidity results.

We would like to point out that certain cotangent disk bundles, such as those of tori $\T^n$, admit symplectic embeddings into closed symplectic manifolds whose Hamiltonian group carries a genuine quasi-morphism, which can be pulled back to yield quasi-morphisms on the Hamiltonian group of these disk bundles. It is an intriguing question whether this pull-back coincides with the restriction of $\mu_0$. A partial result in this direction is presented in \cite{Monzner_Zapolsky_comparison}. In particular, it is unclear whether the quasi-morphism on the disk cotangent bundle of a torus is invariant under coverings, like $\mu_0$ (see proposition \ref{prop_mu_p_invt_under_coverings}).

Next, we mention that the construction of Lagrangian spectral invariants on the Hamiltonian group can be performed for any symplectically aspherical Lagrangian in a completely analogous matter. One only needs to work with the space of paths with endpoints on the Lagrangian which represent a trivial element in the relative $\pi_1$, and consistently introduce spanning half-disks. The rest of the theory comes though. It is appropriate to mention that a related approach, although in a different context, was pursued by R\'emi Leclercq in \cite{Leclercq_Lagrangian_spectral_invariants}. He defined invariants of Lagrangian submanifolds instead of Hamiltonian diffeomorphisms, but in fact his construction allows for a generalization of the results presented here to the case of symplectically aspherical Lagrangians. We have not done this in detail, but it is likely that this would yield results analogous to those listed above, namely, applications to the fragmentation norm, to Hofer and spectral geometry on the Hamiltonian group of the symplectic manifold in which the chosen Lagrangian is contained, to restrictions on Poisson brackets, and to symplectic rigidity.\footnote{There is in general no analogue of Mather's alpha function for aspherical Lagrangians; the present approach would yield something that may be considered a generalization to such manifolds. In particular, spectral invariants, and consequently partial quasi-morphisms can be defined for Hamiltonians with complete flow, and this may be considered as a generalized alpha function (at zero).}

Lastly, we mention a conjecture due to Viterbo concerning a certain bound on Lagrangian spectral invariants, namely, it states that there is a constant $\kappa$ such that if $\phi \in \cG$ is generated by a Hamiltonian whose support is contained in the unit disk bundle, then $\ell_+(\phi)-\ell_-(\phi) \leq \kappa$ (here $\ell_\pm$ are the Lagrangian spectral invariants introduced in section \ref{section_spectral_invariants} below). If this conjecture is true, the triangle inequality and Poincar\'e duality will immediately imply that $\ell_+$ and $\mu_0$ are quasi-morphisms when restricted to the subgroup of $\cG$ generated by Hamiltonians with support inside the unit disk bundle. This will have applications to second bounded cohomology of this subgroup, to asymptotics of the Hamilton-Jacobi equation, more restrictions on Poisson brackets, and more.

\begin{acknow}We would like to thank Fr\'ed\'eric Bourgeois, Michael Entov, Albert Fathi, Vincent Humili\`ere, Joe Johns, Sergei Lanzat, Slava Matveyev, Maxim Maydanskiy, Marco Mazzucchelli, Dusa McDuff, Fabien Ng\^o, Andreas Ott, Sheila Sandon, and Matthias Schwarz for stimulating discussions, and Leonid Polterovich and Karl Friedrich Siburg for reading a preliminary version of the paper and making valuable suggestions.

AM is partially supported by the German National Academic Foundation. NV is partially supported by the ANR grant ``Floer Power'', ANR-08-BLAN-0291-03/04. FZ thanks the Max Planck Institute for Mathematics in the Sciences, Leipzig, where part of this work was carried out, for hospitality and an excellent research atmosphere. AM and FZ profited from their visit to the University of Chicago, and wish to thank Leonid Polterovich for the invitation. We would like to express our collective gratitude to the organizers of Edi-Fest at ETH Z\"urich, where the idea to write a joint paper was born.

And finally, we wish to acknowledge our intellectual debt to the fascinating paper \cite{Viterbo_homogenization}. In fact, the present work grew out of an attempt to understand it. We are grateful to Claude Viterbo for explaining to us some of its more difficult parts.
\end{acknow}

\section{Spectral invariants for Hamiltonian diffeomorphisms}\label{section_spectral_invariants}

In this section we present the construction and properties of Lagrangian and Hamiltonian spectral invariants on the group $\cG$. Subsection \ref{section_Lagr_sp_invts_FH} contains the construction and properties of Lagrangian spectral invariants arising in Floer homology of the zero section $N\subset T^*N$. Subsection \ref{section_Ham_sp_invts} describes Hamiltonian spectral invariants, subsection \ref{section_comparison_Ham_Lagr_sp_invts} compares them to the Lagrangian invariants. In subsection \ref{section_Lagr_sp_invts_gfs} we briefly review Lagrangian spectral invariants coming from generating functions, and their comparison to the Floer-homological ones. Finally, subsection \ref{section_summary_sp_invts} summarizes the various properties of the spectral invariants.

Fix a closed connected manifold $N$. All homology is with $\Z_2$ coefficients, and all moduli spaces are counted modulo $2$. We identify $N$ with the zero section in $T^*N$ via the embedding $N \to T^*N$.

All the material in this section is known and more or less standard, with the exception of the sharp triangle inequality for Lagrangian spectral invariants, proposition \ref{prop_triangle_ineq}, its consequence, the independence of spectral invariants of isotopy, lemma \ref{lemma_indep_of_isotopy}, and the comparison of Lagrangian and Hamiltonian spectral invariants, proposition \ref{prop_comparison_Ham_Lagr_sp_invts}. The exposition is terse, but on the other hand it is extensive enough so as to provide sufficient background both for the sake of proof of the new results, and for the reader who is familiar with Floer homology, but not with spectral invariants.

\subsection{Lagrangian spectral invariants from Floer homology}\label{section_Lagr_sp_invts_FH}

Here we define Lagrangian spectral invariants for Hamiltonian diffeomorphisms via Lagrangian Floer homology and prove some of their properties. The general reference we use is Oh's works \cite{Oh_action_I}, \cite{Oh_action_II}. Whatever statements we make without proof or reference can be found there. We would like to point out that our sign conventions are different from those of Oh. The effect of this difference is that our invariants are ``dual'' to his. This is discussed in subsection \ref{section_sign_conventions}.

The setup is as follows. Let $H \in C^\infty_c([0,1]\times T^*N)$. We define the action functional $\cA_H$ on the space of paths
$$\Omega = \{\gamma\fc[0,1]\to T^*N\,|\,\gamma(0)\in N\}$$
in $T^*N$ by
$$\cA_H(\gamma) = \int_0^1H_t(\gamma(t))\,dt - \int\gamma^*\lambda\,.$$
Let $M \subset N$ be a closed connected submanifold. Consider the path space
$$\Omega(M) = \{\gamma\in\Omega\,|\, \gamma(1)\in\nu^*M\}\,,$$
where $\nu^*M \subset T^*N$ is the conormal bundle of $M$ in $N$. We let $\cA_{H:M}$ be the restriction of $\cA_H$ to $\Omega(M)$. The set $\Crit(H:M) = \Crit \cA_{H:M}$ of critical points of $\cA_H|_{\Omega(M)}$ is precisely the set of solutions $\gamma$ of the Hamiltonian equation of motion, $\dot\gamma=X_H(\gamma)$, with boundary conditions dictated by $\Omega(M)$. The map $\Crit(H:M) \to \phi_H(N)\cap\nu^*M$ given by $\gamma\mapsto\gamma(1)$ is a bijection. We let the action spectrum of $H$ relative to $M$ be the set
$$\Spec(H:M)=\{\cA_H(\gamma)\,|\,\gamma\in\Crit(H:M)\} \subset \R\,.$$
This is a compact nowhere dense subset, and it only depends on the time-$1$ map $\phi_H$ (see, for instance, subsection \ref{section_triangle_ineq}). Consider the vector space $CF(H:M)$ spanned over $\Z_2$ by the set $\Crit(H:M)$, and for $a \notin \Spec(H:M)$, the subspace $CF^{<a}(H:M) \subset CF(H:M)$ spanned by critical points with action $<a$, and the quotient space $CF^{>a}(H:M):=CF(H:M)/CF^{<a}(H:M)$.

Let $J \fc [0,1] \to \End(TT^*N)$ be a path of almost complex structures, compatible with $\omega$ in the sense that $\omega(\cdot,J_t\cdot)$ is a path of Riemannian metrics on $T^*N$. There is an induced $L^2$-metric on $\Omega(M)$, as follows: for $\xi,\eta \in T_\gamma\Omega(M)$ put $\langle \xi, \eta \rangle = \int_0^1\omega(\xi(t),J_t\eta(t))\,dt$. The gradient of $\cA_H$ relative to this metric reads
$$\nabla_\gamma\cA_H(t) = J_t(\gamma(t))\big(\dot\gamma(t)-X_H(\gamma(t)\big)\,.$$
The corresponding negative gradient equation for $u \fc \R(s) \to \Omega(M)$ is Floer's equation
$$\frac{\partial u}{\partial s} + J_t(u) \left(\frac{\partial u}{\partial t} - X_H(u)\right)=0\,.$$
For $\gamma_\pm \in \Crit (H:M)$ we let $\widehat\cM(\gamma_-,\gamma_+)$ denote the set of solutions $u$ of this equation such that $u(\pm\infty,\cdot)=\gamma_\pm$; this set admits a natural action of $\R$ by translation in the $s$ variable, and we let $\cM(\gamma_-,\gamma_+) = \widehat\cM(\gamma_-,\gamma_+)/\R$ be the quotient if $\gamma_+\neq\gamma_-$ and $\cM(\gamma_-,\gamma_-) = \varnothing$.

\subsubsection{Generic Hamiltonian}\label{section_generic_Hamiltonian}

For a generic choice of $H$ the intersection $\phi_H(N) \cap \nu^*M$ is transverse and so $\Crit(H:M)$ is finite, and the various spaces $CF$ are all finite-dimensional; we also refer to such a Hamiltonian as regular. If in addition $J$ is chosen generically, then for any $\gamma_\pm \in \Crit(H:M)$ the moduli spaces $\widehat\cM(\gamma_-,\gamma_+)$, $\cM(\gamma_-,\gamma_+)$ are finite-dimensional smooth manifolds; we also call such a $J$ regular for $H$. There is an integer-valued index $m_{H:M} \fc \Crit(H:M) \to \Z$ such that $\dim \widehat\cM(\gamma_-,\gamma_+) = m_{H:M}(\gamma_-)-m_{H:M}(\gamma_+)$, the Conley-Zehnder index. There are various conventions in the literature concerning its normalization; we use the following one: let $f_0 \fc M \to \R$ be a Morse function; identify a neighborhood of $M \subset N$ with a disk bundle $\pi_0\fc DM \to M$ in the normal bundle $\nu_NM$, extend $\pi_0^*f_0$ to a smooth function $f$ on $N$, and let $H = \pi^*f$. Then elements of $\Crit (H:M)$ are in $1$-$1$ correspondence with the critical points of $f_0$. We normalize $m_{H:M}$ so that it coincides with the Morse index of $f_0$ under this correspondence. We let $CF_k(H:M)$ denote the subspace of $CF(H:M)$ spanned by elements of index $m_{H:M}=k$.

When $J_t$ coincides, outside a compact subset of $T^*N$, with the almost complex structure induced by the Riemannian metric on the base, the various moduli spaces $\cM$ become compact up to breaking. In particular, if $m_{H:M}(\gamma_-) = m_{H:M}(\gamma_+)+1$, $\dim \cM(\gamma_-,\gamma_+)=0$ and so we can define $\partial \fc CF_k(H:M) \to CF_{k-1}(H:M)$ by the linear extension of
$$\partial \gamma_- = \sum_{m_{H:M}(\gamma_+)=k-1}\#\cM(\gamma_-,\gamma_+)\,\gamma_+\,.$$
We have $\partial^2=0$ and the corresponding Floer homology groups are $HF_*(H:M)$. Since elements of $\cM(\gamma_-,\gamma_+)$ are negative gradient lines of the action functional, it decreases along any such element; therefore $\partial$ induces a differential on the subspace $CF^{<a}_*(H:M)$, as well as on the quotient space $CF^{>a}_*(H:M)$. We let $i^a_* \fc HF^{<a}_*(H:M) \to HF_*(H:M)$ and $j^a_* \fc HF_*(H:M) \to HF^{>a}_*(H:M)$ be the induced maps on homology.

The various groups $HF$, as well as the morphisms $i_*^a,j_*^a$, are independent of $J$, which is why we suppressed it from the notation. Moreover, if $K$ is another Hamiltonian, there is a canonical continuation isomorphism $HF_*(H:M) \simeq HF_*(K:M)$. When $H = \pi^*f$ for $f$ a function on $N$ constructed as in the first paragraph of this subsection, the Floer complex of $H$ degenerates into the Morse complex of $f_0$, including grading, which shows that, for any $H$, $HF_*(H:M)$ is canonically isomorphic to the singular homology $H_*(M)$. Using this identification, we can define, for generic $H$, the Lagrangian spectral invariants $\ell(\alpha,H:M)$ for $\alpha \in H_*(M)$ by
$$\ell(\alpha,H:M):=\inf\{a\,|\,\alpha \in \im i^a_*\}\,.$$
These have the following properties, proved by Oh:
\begin{enumerate}
\item $\ell(\alpha,H:M) \in \Spec(H:M)$, in particular it is a finite number;
\item if $H_k$ is a sequence of regular Hamiltonians which tends to $0$ in the $C^1$-topology, then $\ell(\alpha,H_k:M) \to 0$;
\item $\int_0^1 \min(H_t-K_t)\,dt \leq \ell(\alpha,H:M) - \ell(\alpha,K:M) \leq \int_0^1 \max(H_t-K_t)\,dt$; in particular the spectral invariants are Lipschitz with respect to the $C^0$-norm.
\end{enumerate}
We refer to property (iii) as the \emph{continuity} of the spectral invariants.

Similarly, we can define spectral invariants associated to cohomology classes of $M$. To this end, consider the dual Floer complex $CF^*(H:M) = \Hom(CF_*(H:M),\Z_2) \equiv (CF^*(H:M))^*$. The universal coefficient theorem implies that the cohomology of this cochain complex taken with the dual differential $\partial^*$ is canonically isomorphic to the dual of its homology, that is to $(H_*(M))^*$, which with coefficients in a field is the same as the singular cohomology $H^*(M)$. The dual complex is similarly filtered by the action, that is, it increases along the differential. More precisely, we consider the subcomplex $CF^*_{>a}(H:M)$ generated by orbits of action $>a$ and the quotient complex $CF^*_{<a}(H:M) = CF^*(H:M)/CF^*_{>a}(H:M)$. Here we identify the basis of $CF_*$ with the dual basis of $CF^*$, and as a result we have canonical identifications $CF^*_{>a}(H:M) = (CF_*^{>a}(H:M))^*$ and $CF^*_{<a}(H:M) = (CF_*^{<a}(H:M))^*$, and the same for (co)homology. We let $j^*_a \fc HF^*_{>a}(H:M) \to HF^*(H:M)$ and $i^*_a \fc HF^*(H:M) \to HF^*_{<a}(H:M)$ be the maps induced on cohomology by the inclusion and projection maps. We then obtain that the short exact sequence of cochain complexes
$$0\to CF^*_{>a}(H:M) \to CF^*(H:M) \to CF^*_{<a}(H:M)\to 0$$
is dual to the short exact sequence of chain complexes
$$0\to CF_*^{<a}(H:M) \to CF_*(H:M) \to CF_*^{>a}(H:M)\to 0\,,$$
and the induced long exact sequence of cohomologies
$$\dots \to HF^{k-1}_{<a}(H:M) \to HF^k_{>a}(H:M) \xrightarrow{j^k_a} HF^k(H:M) \xrightarrow{i^k_a}HF^k_{<a}(H:M)\to \dots$$
is dual to the long exact sequence of homologies
$$\dots \to HF_{k+1}^{>a}(H:M) \to HF_k^{<a}(H:M) \xrightarrow{i_k^a} HF_k(H:M) \xrightarrow{j_k^a}HF_k^{>a}(H:M)\to \dots\,.$$

The spectral invariant corresponding to $v \in H^*(M)$ is
$$\ell(v,H:M) = \sup\{a\,|\,i^*_a(v) = 0\}\,.$$

\subsubsection{Arbitrary Hamiltonian and the action homomorphism}\label{section_arbitrary_Hamiltonian}

If $H$ is an arbitrary compactly supported Hamiltonian, it can be approximated by regular (that is, generic) Hamiltonians $H_k$, in the $C^\infty$ sense; it follows from the continuity of spectral invariants that $\ell(\alpha,H_k:M)$ is a convergent sequence and that its limit only depends on $H$. Thus spectral invariants can be uniquely extended to the set of all Hamiltonians. It can be proved that these extended invariants satisfy the spectrality axiom (see, for instance \cite{Oh_spectral_invariants}; in our case it is even easier since one does not have to keep track of spanning disks), that is
$$\ell(\alpha,H:M) \in \Spec(H:M)\,.$$
Of course, the extended invariants are also continuous in the sense of property (iii) above, and so they are Lipschitz with respect to the $C^0$-norm.

Let us prove proposition \ref{prop_action_homomorphism} which states that there is a natural homomorphism on the subgroup $\cG_0 \subset \cG$ which consists of Hamiltonian diffeomorpisms fixing the zero section as a set.
\begin{prf}[of proposition \ref{prop_action_homomorphism}]The homomorphism $\cA \fc \cG_0 \to \R$ is defined as follows. Let $H$ be a time-dependent Hamiltonian generating $\phi \in \cG_0$. Pick $q \in N$, put $\gamma_q(t)=\phi_H^t(q)$ and define
$$\cA(\phi)=\cA_H(\gamma_q)\,.$$
Let us first see that the above action does not depend on the choice of the point $q$. Indeed, let $N \to \Omega(N)$ be the map $x \mapsto \gamma_x$ where $\gamma_x(t) = \phi_H^t(x)$. Then it is a smooth embedding and has as its image the set of critical points of $\cA_H$. Since any function attains the same value on a connected submanifold which consists solely of critical points, we see that $\cA_H(\gamma_x)$ is independent of $x$. Thus the isotopy $\phi_H^t$ has as its spectrum only one point. Remark \ref{rem_indep_action_spectrum_of_iso} shows that the action spectrum is independent of the isotopy representing a given element of $\cG$ and therefore $\cA$ is well-defined. It is a homomorphism because action is additive under concatenations. \qed
\end{prf}

As a consequence of spectrality, we have the following observation, which turns out to be crucial for many applications of Lagrangian spectral invariants:
\begin{lemma}\label{lemma_computing_Lagr_sp_invts}The restriction of any spectral invariant $\ell(\alpha,\cdot:N)$ to the group $\cG_0$ coincides with the action homomorphism. It follows that if $H \in C^\infty_c([0,1]\times T^*N)$ satisfies $H|_N \geq c$ (respectively $H|_N\leq c$) for some $c \in \R$, then $\ell(\alpha,H:N) \geq c$ (respectively $\ell(\alpha,H:N) \leq c$), for any $\alpha \neq 0$. In particular, if $H|_N = c$, then $\ell(\alpha,H:N)=c$.
\end{lemma}

\begin{prf}For the first assertion let $H$ be a Hamiltonian generating an element $\phi \in \cG_0$. Spectrality implies that $\ell(\alpha,H:N)$ equals the action of an orbit of the flow of $H$. Proposition \ref{prop_action_homomorphism} shows that this action equals $\cA(\phi)$. This proves the first assertion, and in particular shows that $\ell(\alpha,H:N)$ only depends on $\phi$.

Assume now that $H|_N = c$. The zero section being Lagrangian, the flow of $H$ preserves it, since $H|_N$ is constant. The action of any orbit equals $c$, thus the proof in this particular case is done.

Now if $H|_N \geq c$, we can find another time-dependent Hamiltonian $K$ with compact support which satisfies $H \geq K$ and $K|_N = c$. The claim then follows from the particular case we just considered and the continuity of spectral invariants. The other inequality is proved similarly. \qed
\end{prf}

\begin{rem}In what follows we will need from time to time to use Hamiltonians defined on $[0,\tau]\times T^*N$ with $\tau$ different from $1$. All the preceding constructions are modified in the obvious way, for example, the action functional is now defined on paths $\gamma\fc[0,\tau]\to T^*N$ by $\cA(\gamma)=\int_0^\tau H_t(\gamma(t))\,dt-\int\gamma^*\lambda$, and so on. We will not mention this modification explicitly, and the context will always make clear the domain of definition of Hamiltonians, paths, and action functionals.
\end{rem}

\subsubsection{Poincar\'e duality}\label{section_Poincare_duality}

In this subsection $M=N$. Let $H$ be regular, that is $\phi_H(N)$ intersects $N$ transversely. By standard duality considerations (see \cite{Schwarz}, for example) we obtain
$$\ell(\text{pt},H) = \ell(1,H) \quad\text{and}\quad \ell([N],H) = \ell(\mu_N,H)\,,$$
where $\text{pt}\in H_0(N)$, $[N] \in H_n(N)$, $1 \in H^0(N)$, $\mu_N \in H^n(N)$ are the generators.

Consider the Hamiltonian $\ol H$ defined by $\ol H(t,x) = -H(1-t,x)$. It generates the isotopy obtained from the one generated by $H$ by retracing it backward, that is
$$\phi^t_{\ol H} = \phi^{1-t}_H\phi^{-1}_H\,.$$
The involution $\Omega(N) \to \Omega(N)$, $\gamma \mapsto \ol\gamma=\gamma(1-\cdot)$ establishes a bijection between the sets of critical points of $\cA_H$ and $\cA_{\ol H}$. Moreover, if $J$ is a compatible almost complex structure, regular for $H$, $\ol J(t,\cdot) = J(1-t,\cdot)$ is compatible and regular for $\ol H$, and there is a natural identification of the moduli spaces $\cM(\gamma_-,\gamma_+,H,J)$ and $\cM(\ol\gamma_+,\ol\gamma_-,\ol H,\ol J)$, given by $u \mapsto \ol u$, $\ol u(s,t) = u(-s,1-t)$. Moreover, we have $m_{H:N}(\gamma) = n-m_{\ol H:N}(\ol\gamma)$. It follows that there is a canonical isomorphism
$$CF^*(H:N) = CF_{n-*}(\ol H:N)\,,$$
with the filtrations reversed, that is
$$CF^*_{\cA_H > a}(H:N) = CF_{n-*}^{\cA_{\ol H} < -a}(\ol H:N)\,,$$
for every $a \notin \Spec(H:N)$, since $\cA_H(\gamma) = -\cA_{\ol H}(\ol\gamma)$. We can conclude that
$$\ell(\text{pt},H) = \ell(1,H) = -\ell([N],\ol H)\,.$$
This continues to hold with $H$ replaced by an arbitrary smooth Hamiltonian, due to continuity of spectral invariants. Similarly, we have
$$\ell([N],H) = \ell(\mu_N,H) = -\ell(\text{pt},\ol H)\,.$$

In fact, one can prove, using an argument similar to (and actually, a little simpler than) that of \cite{EP_Calabi_qm}, that the following more general version of Poicar\'e duality holds:
$$\ell(\alpha,H) = -\inf\{\ell(u,\ol H)\,|\,u\in H^*(N)\text{ and }u(\alpha)\neq 0\}\,,$$
or using homology only,
$$\ell(\alpha,H) = -\inf\{\ell(\beta,\ol H)\,|\,\beta\in H_*(N)\text{ and }\alpha\cap\beta\neq 0\}\,.$$
We will not need this more general version, however.

\subsubsection{Triangle inequality and independence of isotopy}\label{section_triangle_ineq}

\begin{notation}\label{notation_ell_pm}In case $M=N$, we denote the corresponding spectral invariants via $\ell(\alpha,H)$ and $\ell(v,H)$. Also, we set $\ell_+=\ell([N],\cdot)$ and $\ell_-=\ell(\text{pt},\cdot)$.
\end{notation}

For the rest of this subsection we assume $M=N$.

Given two functions $H,H'\fc[0,1]\times T^*N \to \R$ such that $H(1,\cdot)=H'(0,\cdot)$, we define their concatenation $H\sharp H' \fc [0,2]\times T^*N$, via
$$H \sharp H'(t,x) = \left\{\begin{array}{ll}H(t,x)\,,&\text{if }t\leq 1\\ H'(t-1,x)\,,&\text{if }t\geq 1\end{array}\right.\,.$$
If $H,H'$ are smooth and $H(1,\cdot)=H'(0,\cdot)$ with all the time derivatives, $H\sharp H'$ is smooth as well.

The first result of this subsection reads
\begin{prop}\label{prop_triangle_ineq}Let $H,H' \in C^\infty_c([0,1]\times T^*N)$ be such that $H(1,\cdot)=H'(0,\cdot)$ with all the time derivatives. Then
$$\ell(\alpha \cap\beta,H\sharp H') \leq \ell(\alpha,H)+\ell(\beta,H')$$
for all $\alpha, \beta \in H_*(N)$ with $\alpha \cap \beta \neq 0$.
\end{prop}
\noindent Here $\cap \fc H_j(N) \times H_k(N) \to H_{j+k-n}(N)$ is the intersection product in homology.

\begin{rem}\label{rem_smoothing}
There is a procedure (see \cite{Polterovich_geom_grp_sympl_diffeo}, for instance) which allows to replace any given time-dependent Hamiltonian with one which vanishes for values of time close to $0$ and $1$, which we call smoothing. This procedure leaves intact all the spectral invariants of the Hamiltonian. Also, the concatenation of any two smoothed Hamiltonians is again smooth. This works as follows.

\begin{enumerate}
\item Consider $H(t,x)$, a time-dependent Hamiltonian on $T^*N$ with compact support. Let $f \fc [0,1]\to[0,1]$ be a smooth function with $f' \geq 0$ everywhere and $f(0)=0$, $f(1)=1$. Let $H^f(t,x)=f'(t)H(f(t),x)$. This is also a smooth Hamiltonian with compact support. Its flows satisfies $\phi_{H^f}^t = \phi_H^{f(t)}$. Thus there is a bijection between the sets of solutions of the corresponding Hamiltonian ODEs with boundary conditions on the zero section, given by $\Crit(H:N) \to \Crit(H^f:N)$, $\gamma \mapsto \gamma^f$, $\gamma^f(t) = \gamma(f(t))$. This bijection preserves the corresponding actions: $\cA_H(\gamma) = \cA_{H^f}(\gamma^f)$. Therefore, if $f_\tau \fc [0,1] \to [0,1]$, $\tau\in[0,1]$, is a continuous family of smooth functions with $f_0=\id_{[0,1]}$, $f_1=f$ and $f_\tau(0) = 0$, $f_\tau(1)=1$, $f_\tau' \geq 0$, then the action spectrum $\Spec(H^{f_\tau}:N)$ is independent of $\tau$, and consequently, by spectrality, so is any spectral invariant.
\item Now let $f$ satisfy the additional requirement that $f(t)=0$ for $t$ near $0$ and $f(t) = 1$ for $t$ near $1$. Consider another time-dependent Hamiltonian $K$ and another function $g$ with the same properties as $f$. The concatenation $H^f \sharp K^g$ is then smooth, and its spectral invariants are independent of the functions $f,g$ used for smoothing; moreover, if the concatenation $H\sharp K$ is smooth, then $H\sharp K$ and $H^f\sharp K^g$ have the same spectral invariants as well. If $H$ is a regular Hamiltonian, then so is $H^f$. If $J$ is an almost complex structure regular for $H$, then $J^f = J(f'(\cdot),\cdot)$ is for $H^f$, with an obvious identification between the various moduli spaces relative to $H,J$ and $H^f,J^f$.
\end{enumerate}
\end{rem}

\begin{prf}[of proposition \ref{prop_triangle_ineq}]
The above remark, together with the continuity of spectral invariants, shows that it suffices to prove the statement for $H,H'$ regular and smoothed, that is, $H=H'=0$ for times $t$ near $0,1$.

Let $\ve > 0$. Consider the concatenation $H''_0 = H \sharp H'$. It may not be regular any more, so we perturb it to a regular Hamiltonian $H''$ such that $\|H''-H''_0\|_{C_0} < \ve$. Moreover, we choose an additional smooth function $K \fc \R \times [0,2] \times T^*N \to \R$ such that $K(s,t, \cdot) = H(t,\cdot)$ for $s \leq 1$ and $t\in[0,1]$, $K(s,t,\cdot)=H'(t-1,\cdot)$ for $s\leq 1$ and $t\in[1,2]$, $K(s,t,\cdot)=H''(t,\cdot)$ for $s\geq2$ and all $t$ and for $s\in[1,2]$ we have $\big|\frac{\partial K} {\partial s}\big| < \ve$ for all $t$.

Fix a $t$-dependent almost complex structure $J$, defined for $t \in [0,2]$, which coincides with the metric almost complex structure outside a compact. For $\gamma,\gamma',\gamma''$ critical points of $\cA_H$, $\cA_{H'}$, $\cA_{H''}$ respectively, we consider the moduli space $\cM(\gamma,\gamma';\gamma'')$ of maps $u \fc \Upsilon \to T^*N$, where $\Upsilon$ is the strip with a slit\footnote{This $\Upsilon$ is a Riemann surface with boundary which is conformally equivalent to a closed disk with three boundary punctures; we put on it the conformal coordinates coming from the identification of its interior with the domain $\R\times(0,2)-(-\infty,0]\times\{1\}\subset \R^2=\C$. The conformal coordinate near the point $(0,1)$ is given by the square root.} appearing in \cite{Abb_Schwarz}, with coordinates $(s,t)$, where $t \in [0,2]$, satisfying
$$\frac{\partial u}{\partial s}(s,t) + J_t(u)\left(\frac{\partial u}{\partial t}(s,t) - X_K(s,t)\right)=0\,,$$
subject to the boundary conditions $u(\partial \Upsilon) \subset N$ and to the asymptotic conditions $u(-\infty,\cdot)=\gamma$, $u(-\infty,\cdot-1)=\gamma'$, $u(\infty,\cdot) = \gamma''$. For a generic choice of $J$, $\cM(\gamma,\gamma';\gamma'')$ is a smooth manifold of dimension $m_{H:N}(\gamma)+m_{H':N}(\gamma')-m_{H'':N}(\gamma'') - n$, compact in dimension $0$. This allows to define a bilinear map
$$CF_j(H:N)\times CF_k(H':N) \to CF_{j+k-n}(H'':N)$$
by the linear extension of
$$(\gamma,\gamma')\mapsto \sum_{\gamma''}\#\cM(\gamma,\gamma';\gamma'')\,\gamma''\,.$$
Examining the boundary of the compactification of the $1$-dimensional such moduli spaces, we see that this bilinear map is in fact a chain map, hence descends to homology,
$$HF_j(H:N)\times HF_k(H':N) \to HF_{j+k-n}(H'':N)\,.$$
We claim that, under the natural identifications $HF_*=H_*(N)$, this map corresponds to the intersection product. Indeed, Oh proved that a different version of this $\Upsilon$-product corresponds to the cup product in singular cohomology. In his version the Hamiltonian $K$ on the strip with a slit vanishes for $s$ near $0$. It can be seen that if we use such a Hamiltonian in the definition of our moduli space, we will obtain the same map on homology. Indeed, one can define the corresponding moduli space of paths of solutions to the above equation where the Hamiltonian depends on the variable of the path, say $K^\tau$. Examining the boundary of the $1$-dimensional such moduli spaces, one can see that counting the $0$-dimensional moduli spaces amounts to a chain homotopy between the chain maps constructed from Hamiltonians $K^0$ and $K^1$, which implies that they define the same map in homology. Thus it is immaterial whether to use our Hamiltonian $K$, ``glued'' from $H,H',H''$, or Oh's Hamiltonian which vanishes for $s$ near $0$. Now, Oh's sign conventions make his Floer homologies isomorphic to $H^*(N)$ (see subsection \ref{section_sign_conventions}). Passing to our sign conventions amounts to applying the Poincar\'e duality in each variable, which transforms the cup product on cohomology into the intersection product on homology.

Now, a computation  shows (compare with \cite{Abb_Schwarz}) that if $u \in \cM(\gamma,\gamma';\gamma'')$, then
$$\cA_H(\gamma)+\cA_{H'}(\gamma') - \cA_{H''}(\gamma'') \geq E(u) - \ve\,,$$
where $E(u) \geq 0$ is the energy of $u$. It follows that the above chain map restricts to a map on filtered subcomplexes:
$$CF_j^{<a}(H:N)\times CF_k^{<b}(H':N) \to CF_{j+k-n}^{a+b+\ve'}(H'':N)$$
for any $a,b,\ve'$ such that $a\notin\Spec(H:N)$, $b \notin \Spec(H':N)$, $\ve' > \ve$, and $a+b+\ve' \notin \Spec(H'':N)$. This implies that
$$\ell(\alpha \cap\beta,H'') \leq \ell(\alpha,H)+\ell(\beta,H')+\ve\,.$$
Since $H''$ was chosen $\ve$-close to the concatenation $H\sharp H'$, passing to the limit as $\ve \to 0$, we obtain the desired triangle inequality
$$\ell(\alpha \cap\beta,H\sharp H') \leq \ell(\alpha,H)+\ell(\beta,H')\,.\qed$$
\end{prf}

As a consequence, we have
\begin{lemma}\label{lemma_indep_of_isotopy}Let $H,H'\in C^\infty_c([0,1]\times T^*N)$ have the same time-$1$ map, $\phi_H = \phi_{H'}$. Then the spectral invariants of $H,H'$ coincide.
\end{lemma}

\begin{prf}
First, let $G\in C^\infty_c([0,1]\times T^*N)$ be a Hamiltonian generating a loop, that is $\phi_G = \id$. We claim that its spectral invariants all vanish. First, observe that we may replace $G$ by a smoothed version, without altering the spectral invariants, and such that $\phi_G$ is still the identity map. Note that $\ell(\alpha,G)$ is, by spectrality, the action of a Hamiltonian arc $\gamma \in \Omega(N)$. Since $G$ generates a loop and is smoothed, this arc is in fact a smooth closed orbit. A standard computation shows (see \cite{Schwarz}) that the actions $\cA_G(\gamma_x)$ are all the same, where $\gamma_x(t)=\phi_G^tx$. It follows that they are all zero, because we can take $x$ to be outside the support of $G$. Thus $\cA_G(\gamma) = \cA_G(\gamma_{\gamma(0)})=0$, as claimed.

Let us now prove that $\ell(\alpha,H) = \ell(\alpha,H')$. Again, assume that $H,H'$ are smoothed by the above procedure so that both equal $0$ near $t=0,1$. Suppose for a moment that we can show the following equality:
$$\ell(\alpha,H) = \ell(\alpha,H\sharp\ol {H'}\sharp H')\,.$$
Then we have
$$\ell(\alpha,H) = \ell(\alpha\cap[N],H\sharp\ol {H'} \sharp H') \leq \ell(\alpha,H') + \ell([N],H\sharp\ol {H'})\,.$$
Since $H\sharp \ol {H'}$ generates a loop, its spectral invariants vanish and we obtain
$$\ell(\alpha,H) \leq \ell(\alpha, H')\,,$$
and the reverse inequality follows by exchanging $H$ and $H'$.

To prove that
$$\ell(\alpha,H) = \ell(\alpha,H\sharp\ol {H'} \sharp H')\,,$$
we proceed as follows. Since we smoothed $H'$, it is true that $H'(t,\cdot)\equiv0$ for $t\in[0,\delta]\cup[1-\delta,1]$ for some $\delta > 0$. Let $f \fc [0,1] \to \R$ be a smooth function such that $f(t)=t$ for $t \in [0,1-\delta/2]$, $f'\geq 0$ everywhere and $f(t)$ is constant on $[1-\delta/4,1]$. Define $f^\tau \fc [0,\tau] \to \R$ by $f^\tau(t) = f(t+\tau)-\tau$, for $\tau \geq \delta$ and $f^\tau \equiv 0$ for $\tau < \delta$. Put $K^\tau(t,x)=(f^\tau)'(t)H'(f^\tau(t),x)$ for $t\in[0,\tau]$ and $\ol{K^\tau}(t,x)=K^\tau(\tau-t,x)$. It is easy to see that for all $\tau\in[0,1]$ the Hamiltonians $K^\tau,\ol{K^\tau}$ are smooth. Now let $H^\tau$ be the concatenation of $H$, then $K^\tau$ running in time $\tau$ and then $\ol{K^\tau}$ running in time $\tau$. An immediate computation shows that $\Spec(H^\tau:N)$ is independent of $\tau$ and that $H^0=H$ and $H^1=H\sharp\ol{H'}\sharp H'$. The assertion now follows from spectrality. \qed
\end{prf}

Henceforth we denote by $\ell(\alpha,\phi)$ the value $\ell(\alpha,H)$ for any $H$ generating $\phi$.

\begin{rem}\label{rem_indep_action_spectrum_of_iso}The fact that the spectrum $\Spec(H:M)$ only depends on the time-$1$ map of $H$ can be proved in a similar, though much more elementary, way, since no triangle inequality is needed. To wit, as we mentioned in the beginning of the proof of lemma \ref{lemma_indep_of_isotopy}, for a Hamiltonian generating a loop the action of any Hamiltonian orbit vanishes. Now let $H,H'$ have the same time-$1$ map, and be smoothed, without loss of generality. Let $z \in T^*N$ and $\gamma(t) = \phi^t_{H}(z)$, $\gamma'(t)=\phi_{H'}^t(z)$, and let $\gamma''$ be the concatenation of $\gamma$ and the reversal of $\gamma'$. Then
$$\cA_H(\gamma) - \cA_{H'}(\gamma') = \cA_{H\sharp \ol{H'}}(\gamma'') = 0\,,$$
since $\gamma''$ is an orbit of $H\sharp \ol{H'}$, which generates a loop.
\end{rem}

\subsubsection{Hamiltonians with complete flow}\label{section_geom_bdd_Hamiltonians}

Here we describe how to define the various spectral invariants $\ell(\alpha,\cdot:M)$ for a Hamiltonian having complete flow.

\begin{lemma}Let $H,H'$ be two time-dependent Hamiltonians with compact support. Assume that there are two open subsets $U\subset V \subset T^*N$ such that $N \subset U$, $\phi^t_H(U),\phi^t_{H'}(U) \subset V$ for all $t\in[0,1]$ and $H|_{[0,1]\times V}=H'|_{[0,1]\times V}$. Then the spectral invariants of $\phi_H$ and $\phi_{H'}$ coincide.
\end{lemma}
\begin{prf}This follows from the fact that $H$ can be continuously deformed into $H'$ such that the action spectrum stays intact during the deformation. More precisely, let $H^\tau = \tau H' + (1-\tau)H$. Then $H^\tau$ is a smooth Hamiltonian whose flow sends $U$ into $V$ for all times and which coincides with $H$ and $H'$ when restricted to $V$. It follows that $H^\tau$ has the same set of Hamiltonian orbits in $\Omega(M)$ regardless of $\tau$ and those have actions independent of $\tau$. The claim follows. \qed
\end{prf}

If $H$ has complete flow, there is $R>0$ such that $\phi_H^t(N) \subset T^*_{<R}N$ for all $t\in[0,1]$. Any two compactly supported cutoffs $H',H''$ of $H$ outside $T^*_{<R}N$ satisfy the assumptions of the lemma and so have identical spectral invariants; we declare the common value $\ell(\alpha,H':M)=\ell(\alpha,H'':M)$ to be the spectral invariant $\ell(\alpha,H:M)$. Note that these extended spectral invariants share the properties of the usual ones, that is, spectrality, continuity, the triangle inequality, independence of isotopy, and, what is also important in applications, the product formula below, for which, incidentally, Hamiltonians with complete flow provide natural subjects.

For future use, we formulate
\begin{lemma}\label{lemma_sp_invts_geom_bdd}Let $H$ be a time-dependent Hamiltonian with complete flow and assume that this flow keeps the zero section inside an open set $U$ for all times. Then if $G$ is any cutoff of $H$ outside $U$, we have for any $t$, $\ell(\alpha,\phi_G^t:M)=\ell(\alpha,\phi_H^t:M)$.\qed
\end{lemma}

\begin{rem}A word of warning is in order. It is \emph{not} true that one can consistently define the Floer complex for a Hamiltonian with complete flow, since moduli spaces of Floer trajectories may fail to be compact without additional assumptions on the behavior of the Hamiltonian at infinity, such as quadratic growth or similar. It is also not true that the Floer complexes of two cutoffs are isomorphic. What is true, and this is what makes the whole theory work, is that the Floer complex of \emph{any} cutoff is well-defined, and that the complexes of different cutoffs are related by canonical chain maps (continuation morphisms) which descend to level-preserving isomorphisms on homology.
\end{rem}

\subsubsection{The product formula}

In this subsection we prove the product formula for spectral invariants, which turns out to be important for applications to symplectic rigidity. Recall the definition of the direct sum of two time-dependent Hamiltonians, subsection \ref{section_main_result}.
\begin{thm}\label{thm_product_formula_sp_invts}Let $H$, $H'$ be time-dependent Hamiltonians with complete flows on $T^*N$, $T^*N'$, respectively. Then, for any $\alpha \in H_*(N)-\{0\}$ and $\alpha'\in H_*(N')-\{0\}$ we have
$$\ell(\alpha\otimes\alpha',H\oplus H') = \ell(\alpha,H) + \ell(\alpha',H')\,,$$
where $\alpha\otimes\alpha'\in H_*(N)\otimes H_*(N')=H_*(N\times N')$.
\end{thm}

Before passing to the proof, we need some preparations. By definition, a filtered graded chain complex is a quadruple $\cV=(V,\vec v,\cA,\partial)$, where $\vec v = (v_1,\dots,v_k)$ is a graded finite set, $V=\Z_2\otimes \vec v$ is the $\Z_2$-vector space spanned by $\vec v$, $\cA\fc \vec v \to \R$ is a $1$-$1$ function, called the action, and $\partial \fc V \to V$ is a differential, which lowers the grading by $1$, and respects the action filtration, that is it preserves $V^{<a}:=\Z_2\otimes(\vec v \cap \{\cA < a\})\subset V$ for every $a \in \R$. Following the usual procedure, one can define the spectral invariants of $\cV$ relative to homology classes in $H(V,\partial)$, which we denote by $\ell(\alpha,\cV)$ for $\alpha \in H(V,\partial) - \{0\}$. Given two filtered graded chain complexes $\cV=(V,\vec v,\cA,\partial)$ and $\cV'=(V',\vec v', \cA', \partial')$, one can form the product filtered graded chain complex $\cV'':=\cV \otimes \cV' = (V'':=V\otimes V',\vec v'':=\vec v \times \vec v', \cA'':=\cA\oplus\cA',\partial'':=\partial\otimes\id_{V'}+\id_V\otimes\partial')$, where $(\cA\oplus\cA')(v_i,v'_j)=\cA(v_i)+\cA(v_j')$, provided that $\cV$ and $\cV'$ are in general position, meaning that $\cA''$ is still $1$-$1$. Then $H(V'',\partial'')=H(V,\partial)\otimes H(V',\partial')$.

The spectral invariants of filtered graded chain complexes satisfy the following product property:
\begin{lemma}\label{lemma_product_formula_sp_invts_chain_cplxs}Let $\cV,\cV',\cV''$ be as above. Then for $\alpha \in H(V,\partial)-\{0\}$, $\alpha'\in H(V',\partial')-\{0\}$ it is true that
$$\ell(\alpha'',\cV'') = \ell(\alpha,\cV)+\ell(\alpha,\cV')\,,$$
where $\alpha'' = \alpha\otimes\alpha' \in H(V'',\partial'')$. \qed
\end{lemma}
\noindent The proof, though elementary, is somewhat involved, and can be extracted from \cite{EP_rigid_subsets}.

We can now pass to the proof of theorem \ref{thm_product_formula_sp_invts}.
\begin{prf}[of theorem \ref{thm_product_formula_sp_invts}]Given an arbitrary Hamiltonian with complete flow, we can always perturb it (say, in $C^\infty$ topology) to a generic one, meaning that the Floer complex of any cutoff is a filtered graded chain complex in the sense of the discussion above. Moreover, given two such Hamiltonians, we can perturb both of them in such a way that both the perturbations and their direct sum are generic. Since spectral invariants are continuous with respect to $C^0$ norm, it suffices to restrict attention to Hamiltonians $H,H'$ which are generic in this sense, and such that the sum $H\oplus H'$ is generic as well, which is what we choose to do.

Let $G,G'$ be cutoffs of $H,H'$. The direct sum $G\oplus G'$ has complete flow. We choose regular almost complex structures $J,J'$ on $T^*N,T^*N'$, respectively, which coincide, outside a large compact, with the metric almost complex structures. We let $R>0$ be large enough so that $T^*_{<R}(N\times N')$ contains the product $T^*_{<r}N\times T^*_{<r'}N'$, where $r$ is large enough so that $T^*_{<r}N$ contains the images of all the critical points of $\cA_G$, as well as the images of the Floer trajectories between pairs of critical points of $\cA_G$ of index difference $1$, and similarly for $r'$, $T^*N'$ and $G'$. Let $G''$ be a cutoff of $G\oplus G'$ outside $T^*_{<R}(N\times N')$. Then it is also a cutoff of $H\oplus H'$, in particular
$$\ell(\alpha'',G'')=\ell(\alpha'',H\oplus H')\,.$$
Moreover, $G''$ is generic by construction, and $J'':=J\oplus J'$ is a regular almost complex structure. It then follows that the Floer complex of $G''$ relative to $J''$ is a filtered graded chain complex, which is the product of the Floer complexes of $G,G'$. Applying lemma \ref{lemma_product_formula_sp_invts_chain_cplxs}, we see that
$$\ell(\alpha\otimes\alpha',H\oplus H')=\ell(\alpha\otimes\alpha',G'')=\ell(\alpha,G)+\ell(\alpha',G')=\ell(\alpha,H)+\ell(\alpha,H')\,.$$
This concludes the proof. \qed
\end{prf}

\subsection{Hamiltonian spectral invariants}\label{section_Ham_sp_invts}

These were defined for weakly exact symplectic manifolds convex at infinity in \cite{Frauenfelder_Schlenk_Ham_convex}, and in the more general setting of semipositive symplectic manifolds convex at infinity in \cite{Lanzat_thesis}. We only present a sketch of the construction, referring the reader to the aforementioned sources for details.

Standard Floer homology cannot be correctly defined for compactly supported Hamiltonians because they are degenerate. To circumvent this difficulty, one considers Hamiltonians which have support in some fixed cotangent ball bundle and which have a certain prescribed behavior near the boundary.

In more detail, fix $R > 0$ and a smooth function $h\fc(-\ve,\infty) \to \R$, where $\ve >0$,  such that $h(t) = 0$ for $t \geq 0$ and $h'(t) \geq 0$, for $t \leq 0$. Moreover, $h'(t)$ should be small enough so that the flow of $h(\|p\|-R)$ does not have non-constant periodic orbits of period $\leq 1$ for $\|p\|\in(-\ve,0)$. Then for $H_t \in C^\infty_c(T^*N)$ such that $H_t(q,p)=h(\|p\|-R)$ for $\|p\|\geq R-\ve$ the Floer complex $CF(H)$ is well-defined if we take as its generators the $1$-periodic orbits of $H$ inside $T^*_{<R}N$ and all of them are non-degenerate. The boundary operator counts, as usual, Floer cylinders running between pairs of such orbits, and the behavior of $H$ near $\|p\|=R$ guarantees that all of them are contained in $T^*_{\leq R-\ve}N$. The Floer complex is graded by the Conley-Zehnder index\footnote{It is normalized so as to equal the Morse index of critical points of a $C^2$-small Hamiltonian, considered as $1$-periodic orbits.} $m_H$, and the boundary operator lowers it by $1$. Let us denote the homology of $CF_*(H)$ by $HF_*(H;h,R)$. There is a PSS-type isomorphism $HF_*(H;h,R)=H_*(T^*N)$. Since $HF_*(H;h,R)$ is filtered by action, spectral invariants are defined in the standard fashion, namely consider the inclusion morphism $i^a\fc HF_*^{<a}(H;h,R) \to HF_*(H;h,R)$, where $HF_*^{<a}(H;h,R)$ is the homology of the subcomplex $CF_*^{<a}(H) \subset CF_*(H)$ spanned by orbits of action $<a$. Then for $\alpha \in H_*(T^*N)$ we can define
$$c(\alpha,H;h,R) = \inf\{a\,|\,\alpha\in\im i^a\}\,.$$
These spectral invariants satisfy all the standard properties, including Lipschitz continuity, triangle inequality, and spectrality. This implies that if $H$ is an arbitrary compactly supported Hamiltonian, and we $C^0$-approximate it by non-degenerate Hamiltonians $H_k$, $k\in\N$, whose behavior for $\|p\| \in [R-\ve_k,\infty)$ is prescribed by the function $h$ as above, and $\ve_k \to 0$, then the sequence $c(\alpha,H_k;h,R)$ is Cauchy and we declare its limit to be the spectral invariant $c(\alpha,H;h,R)$.

It can be shown, by a standard but a little lengthy argument 
that this spectral invariant is independent of the choices, that is, $h$ and $R$. Moreover, it actually only depends on $\phi_H$, and so we will use the notation $c(\alpha,\phi_H)$ for it. We will only need the invariant $c_-(\phi) = c(\text{pt},\phi)$ and its dual counterpart $c_+(\phi)=-c_-(\phi^{-1})$.

In \cite{Frauenfelder_Schlenk_Ham_convex} it is also shown that if $U \subset T^*N$ is displaceable by $\psi \in \cG$ then for any $\phi \in \cG_U$ it is true that
$$-\Gamma(\psi) \leq c_-(\phi) \leq c_+(\phi) \leq \Gamma(\psi)$$
where $\Gamma(\psi)=c_+(\psi)-c_-(\psi)$ is the spectral norm of $\psi$.

\begin{rem}We would like to point out that Hamiltonian spectral invariants, unlike Lagrangian ones, \emph{are not} defined for Hamiltonians with complete flow. This has to do with the fact that Floer homology constructed from periodic orbits, may be ill-defined for such Hamiltonians. Moreover, there is no consistent way of cutting such Hamiltonians off, like in the Lagrangian case, in order to use the compactly supported theory.
\end{rem}

\subsection{Comparison of Lagrangian and Hamiltonian spectral invariants}\label{section_comparison_Ham_Lagr_sp_invts}

Our goal in this subsection is the following proposition.
\begin{prop}\label{prop_comparison_Ham_Lagr_sp_invts}Let $\phi \in \cG$. Then
$$\ell_-(\phi) \geq c_-(\phi)\,.$$
\end{prop}
\noindent This implies the following chain of inequalities:
$$c_-(\phi) \leq \ell_-(\phi) \leq \ell_+(\phi) \leq c_+(\phi)\,,$$
where the rightmost inequality follows by duality.

\begin{prf}
The proof is essentially contained in \cite{Albers_comparison}, the only point of difference being that there the theory is restricted to closed manifolds. This has to do with compactness of moduli spaces of perturbed pseudo-holomorphic curves. In our case, since the almost complex structure is assumed to coincide, outside a large compact, with the one coming from the auxiliary Riemannian metric, there are no additional compactness issues beyond the closed case, and in fact, since the form is exact, there is no bubbling off of spheres or disks, and the proofs are actually much simpler, and in particular the PSS morphisms constructed in the aforementioned paper are defined for all degrees and are isomorphisms. Therefore we only present a sketch of the argument, emphasizing the essential point of comparison of the spectral invariants.

Albers defines a map $\iota \fc CF_*(H:N) \to CF_*(H)$, as follows. First, one can assume that $H_t$ is time-independent near $t=0,1$, and that the Floer homology for it is defined as in subsection \ref{section_Ham_sp_invts}. Given a Hamiltonian arc $\gamma$ and a periodic orbit $x$ of $H$, consider the moduli space $\cM(\gamma,x)$ consisting of solutions of the Floer equation defined on the Riemann surface $\Upsilon'$, conformal to a closed disk with one boundary and one interior puncture, obtained from the above strip with a slit $\Upsilon$ through identifying the top and the bottom boundary components, that is, solutions to
$$\frac{\partial u}{\partial s} + J_t(u)\left(\frac{\partial u}{\partial t} - X_H\right)=0\,,$$
where the boundary puncture is asymptotic to $\gamma$ and the interior puncture is asymptotic to $x$, while the boundary is mapped to the zero section. What makes this equation well-defined is the presence of global conformal coordinates $(s,t)$ on $\Upsilon'$. He shows that this moduli space is a smooth manifold of dimension $m_{H:N}(\gamma) - m_H(x)$. It is compact in dimension $0$, which follows from the usual convexity considerations. Let $\iota$ be the linear extension of
$$\iota(\gamma)=\sum_{m_{H}(x) = m_{H:N}(\gamma)}\#\cM(\gamma,x)\,x\,.$$
He then shows that this is a chain map, and the canonical identifications $HF_*(H:N) = H_*(N)$ and $HF_*(H) = H_*(T^*N)$ intertwine it with the isomorphism $H_*(N) \to H_*(T^*N)$ induced by the inclusion of the zero section into $T^*N$. We have the sharp action-energy identity for an element $u \in \cM(\gamma,x)$:
$$\cA_H(\gamma) - \cA_H(x) = E(u) \geq 0\,,$$
whence $\iota$ actually maps $CF_*^{<a}(H:N) \to CF_*^{<a}(H)$ for $a\notin\Spec(H)\cup\Spec(H:N)$. It follows that $c_-(H) \leq \ell_-(H)$. Using the continuity of spectral invariants, we conclude that this inequality holds for all smooth Hamiltonians with compact support. \qed
\end{prf}

It follows from the previous subsection that if $U$ is an open subset displaceable by $\psi$, then for any $\phi \in \cG_U$ we have
$$-\Gamma(\psi) \leq \ell_-(\phi) \leq \ell_+(\phi) \leq \Gamma(\psi)\,.$$

\subsection{Lagrangian spectral invariants from generating functions}\label{section_Lagr_sp_invts_gfs}

We also need to use another definition of Lagrangian spectral invariants, namely those coming from generating functions, due to Viterbo. The reason is that we need both definitions in the proof theorem \ref{thm_our_constr_equals_homogenization} which states that the symplectic homogenization is a particular case of our functionals $\mu_a$.

\subsubsection{Definition}\label{section_definition_Lagr_sp_invts_gf}

Let us recall briefly Viterbo's construction of Lagrangian spectral invariants using generating functions \cite{Viterbo_gfqi}. A generating function quadratic at infinity, or gfqi for short, is a function $S \fc N(q) \times E(\xi) \to \R$, where $E$ is a finite-dimensional vector space, such that $\|\partial_\xi S - \partial_\xi B\|_{C^0}$ is bounded, where $B \fc E \to \R$ is a non-degenerate quadratic form. We let $E = E^+ \oplus E^-$ be the splitting into the positive and negative subspaces of $B$.

Consider the relative homology $H_*(\{S<a\},\{S<b\})$. It follows from the definition that for $a$ large enough and $b$ small enough this group is independent of $a,b$ and is canonically isomorphic to $H_*(N) \otimes H_*(E^-,E^--0) \simeq H_{*+d}(N)$, where $d = \dim E^-$ and the last isomorphism (``Thom isomorphism'') is given by tensoring with the generator of $H_d(E^-,E^--0) \simeq \Z_2$. We denote this group by $H_*(S:N)$. There is a natural inclusion morphism $i^b\fc H_*(\{S < b\}) \to H_*(S:N)$. To each $\alpha \in H_*(N)$ we associate the spectral invariant
$$\ell(\alpha,S) = \inf \{b\,|\,\alpha\in\im i^b\}\,.$$
Similarly, if $M \subset N$ is a closed submanifold, we can consider the restriction $S|_{M\times E}$ as a gfqi with base $M$ and define spectral invariants associated to classes in $H_*(M)$.

These invariants are defined also for Lagrangian submanifolds of $T^*N$, as follows. A regular gfqi gives rise to a Lagrangian immersion, see \cite{Viterbo_gfqi}. A Lagrangian submanifold Hamiltonian isotopic to the zero section admits a gfqi unique up to gauge a transformation, stabilization, and the addition of a constant \cite{Viterbo_gfqi}, \cite{Theret_uniqueness_gfqi}. Except the addition of a constant, the elementary operations do not alter the spectral invariants. We can then say that the spectral invariants of this gfqi are attached to the Lagrangian submanifold in question, and they are all defined up to simultaneous addition of a constant.

\subsubsection{Sign conventions}\label{section_sign_conventions}

In what follows we will rely on certain results due to Milinkovi\'c and Oh on the equality of spectral invariants coming from Floer homology and from generating functions \cite{Milinkovic_Oh_FH_as_stable_MH}, \cite{Milinkovic_Oh_gf_versus_action}. Since our sign conventions differ from theirs, it is necessary to relate the two conventions regarding spectral invariants.

Our sign convention follows the philosophy that the Floer theory of the action functional is a perturbation of the Morse theory of a function on a closed manifold, in particular the Hamiltonian enters the action functional with a positive sign.

\begin{notation}In this subsection, as well as in the rest of the paper, we will denote objects defined with the sign conventions of Milinkovi\'c and Oh, with the overline, with the exception of $\ol H$, which we reserve for the ``reversed'' Hamiltonian.
\end{notation}

Namely, assume that $H$ is a compactly supported time-dependent Hamiltonian on $T^*N$. The action functional $\ol\cA_H\fc\Omega(M) \to \R$ is defined as
$$\ol\cA_H (\gamma) = -\cA_H(\gamma) = -\int_0^1H_t(\gamma(t))\,dt+\int\gamma^*\lambda\,.$$
The symplectic form $\ol \omega = -\omega=-d\lambda=-dp\wedge dq$. The Hamiltonian vector field $\ol X_H$ is defined by the equation $\ol\omega(\ol X_H,\cdot)=dH$ and so $\ol X_H = X_H$. In particular, the flows in the two sign conventions coincide.

Assume now that $H$ is regular, that is $\phi_H(N)$ intersects $\nu^*M$ transversely. Then, of course, $\Crit(\ol\cA_H:M) = \Crit(H:M)$, while the action spectrum is flipped: $\Spec(\ol\cA_H:M) = -\Spec(H:M)$. Milinkovi\'c and Oh use the negative gradient flow of $\ol\cA_H$ to produce the Floer equation. In their sign conventions, an almost complex structure $\ol J$ is compatible with $\ol\omega$ if $\ol\omega(\cdot,\ol J\cdot)$ is a Riemannian metric. This is the case if and only if the almost complex structure $J = -\ol J = (\ol J)^{-1}$ is compatible with $\omega$ in our sense. Therefore their negative gradient flow corresponds to our positive gradient flow. It follows that there is a canonical identification of moduli spaces
$$\ol\cM(\gamma_+,\gamma_-) = \cM(\gamma_-,\gamma_+)$$
for $\gamma_\pm \in \Crit(H:M)$. Consequently their Floer boundary operator is the dual of ours. Their convention for the Conley-Zehnder index is $\ol m_{H:M}(\gamma)=\dim M-m_{H:M}(\gamma)$ for $\gamma \in \Crit(H:M)$. Thus their Floer complex
$$(\ol {CF}_*(H:M),\ol\partial_{H:M})$$
is canonically isomorphic to
$$(CF^{\dim M-*}(H:M),(\partial_{H:M})^*)\,,$$
and so the homology they obtain is in fact the singular cohomology $H^{\dim M - *}(M)$, which by Poincare duality is isomorphic to $H_*(M)$. Using this latter identification and the fact that $\ol\partial$ decreases the action $\ol\cA_H$, they define spectral invariants by the usual recipe; we denote them by $\ol\ell(\alpha,H:M)$ for $\alpha \in H_*(M)$.

We will only need the relation between $\ell$ and $\ol\ell$ in the case $M=N$. Then we have, by Poincare duality described in subsection \ref{section_Poincare_duality},
$$(\ol {CF}_*(H:N),\ol\partial_{H:N})=(CF^{n-*}(H:N),(\partial_{H:N})^*)=(CF_*(\ol H:N),\partial_{\ol H:N})\,,$$
and since we have $\ol \cA_H(\gamma)=\cA_{\ol H}(\ol\gamma)$, meaning that the action filtration on $(\ol {CF}_*(H:M),\ol\partial_{H:M})$ by $\ol\cA$ coincides with that on $(CF_*(\ol H:M),\partial_{\ol H:M})$ by $\cA_{\ol H}$, the sought-for relation between the spectral invariants is given by
$$\ol\ell(\text{pt},H)=\ell(\text{pt},\ol H)=-\ell([N],H)\text{\quad and \quad}\ol\ell([N],H)=\ell([N],\ol H)=-\ell(\text{pt},H)\,.$$
Dualizing, spectral invariants corresponding to cohomology classes can be defined and for these we have
$$\ol\ell(1,H)=-\ell(\mu_N,H)\text{\quad and \quad}\ol\ell(\mu_N,H)=-\ell(1,H)\,.$$
Using the succinct notation \ref{notation_ell_pm}, we can write these relations as
$$\ol\ell_\pm(H)=-\ell_\mp(H)\,.$$

\subsubsection{Relation between Floer-homological invariants and invariants from generating functions}\label{section_relation_FH_gf_invts}

In \cite{Milinkovic_Oh_FH_as_stable_MH}, \cite{Milinkovic_Oh_gf_versus_action} Milinkovi\'c and Oh show that the Lagrangian invariants coming from Floer homology and those coming from generating functions coincide provided the generating function is suitably normalized. The normalization is as follows. Recall that if a generating function $W$ defined on the total space of a submersion $\pi \fc E \to N$ generates a Lagrangian \emph{embedding} $L \subset T^*N$, then it induces a function, denoted $W|L$, on the image of the embedding via the formula $W|L := W\circ (i_W)^{-1} \fc L \to \R$, where $i_W\fc\Sigma_W \to T^*N$ is the canonical map from the fiberwise critical locus $\Sigma_W$ of $W$ to $T^*N$. It follows from the fact that the differential $d(W|L)$ coincides with $\lambda|_L$ and that two generating functions of $L$ induce functions on it whose difference is constant. A particularly important case is the action functional $\ol\cA_H$, defined on the path space $\Omega$, where the submersion is given by $\Omega \to N$, $\gamma \mapsto \gamma(1)$. It generates $\phi_H(N)$.

We have the following lemma, whose proof can be extracted from \cite{Milinkovic_Oh_FH_as_stable_MH}, \cite{Milinkovic_Oh_gf_versus_action}:
\begin{lemma}\label{lemma_sp_invts_FH_gfqi_coincide}Let $L \subset T^*N$ be a Lagrangian submanifold, Hamiltonian isotopic to the zero section, and assume that $H_t$ is a compactly supported Hamiltonian such that $\phi_H(N)=L$, and that $S \fc N \times E \to \R$ is a gfqi generating $L$. If the induced functions $\ol\cA_H|L$ and $S|L$ are equal, then for any closed connected submanifold $M \subset N$ and any $\alpha \in H_*(M)$ we have
$$\ol\ell(\alpha,H:M)=\ell(\alpha,S|_{M \times E})\,. \qed$$
\end{lemma}
\noindent In order to see that our normalization condition implies the conclusion of Milinkovi\'c and Oh, one has to translate it into the language of wavefronts, see \cite{Oh_action_I}. It suffices to note that when $S$ is normalized as in the lemma, it is possible to find an interpolation between $S$ and $\cA_H$, in the sense of \cite{Milinkovic_Oh_FH_as_stable_MH}, \cite{Milinkovic_Oh_gf_versus_action}, which has a constant wavefront, and this allows their argument to work.

\subsection{Summary}\label{section_summary_sp_invts}

Here we summarize for further reference the properties of spectral invariants proved above, together with some immediate consequences.

\begin{thm}\label{thm_summary_spectral_invts}Let $N$ be a closed connected manifold. To each $\alpha \in H_*(N) - \{0\}$ we associate a function $\ell(\alpha,\cdot) \fc \cG \to \R$ such that:
\begin{enumerate}
\item $\ell(\alpha,\phi) \in \Spec(\phi:N)$;
\item if $H$, $K$ generate $\phi$, $\psi$, then $\int_0^1 \min(H_t-K_t)\,dt \leq \ell(\alpha,\phi) - \ell(\alpha,\psi) \leq \int_0^1 \max(H_t-K_t)\,dt$;
\item $\ell(\alpha\cap\beta,\phi\psi) \leq \ell(\alpha,\phi)+\ell(\beta,\psi)$; in particular, $\ell_+(\phi\psi) \leq \ell_+(\phi) + \ell_+(\psi)$;
\item $\ell_-(\phi) \leq \ell(\alpha,\phi) \leq \ell_+(\phi)$;
\item $\ell_\pm(\phi) = -\ell_\mp(\phi^{-1})$, and thus $\ell_-(\phi\psi) \geq \ell_-(\phi) + \ell_-(\psi)$;
\item the restriction of $\ell(\alpha,\cdot)$ to $\cG_0$ coincides with the action homomorphism $\cA$; in particular if $H$ generates $\phi$ and $H|_N = c$ (respectively $H|_N \geq c$, $H|_N \leq c$) for some $c \in \R$, then $\ell(\alpha,\phi) = c$ (respectively $\ell(\alpha,\phi) \geq c$, $\ell(\alpha,\phi) \leq c$);
\item if $U \subset T^*N$ is an open subset and $\psi \in \cG$ is such that $\psi(U) \cap U = \varnothing$, then
$$-\Gamma(\psi)\leq \ell_-(\phi) \leq \ell_+(\phi)\leq \Gamma(\psi)$$
for $\phi \in \cG_U$;
\item $|\ell(\alpha,\phi)-\ell(\alpha,\psi\phi\psi^{-1})| \leq \ell_+(\psi)-\ell_-(\psi)$.
\end{enumerate}
\end{thm}

\begin{prf}With the exception of points (iv) and (viii), these statements are proved in the previous subsections. For point (iv) the triangle inequality implies, for example:
$$\ell(\alpha,\phi) = \ell(\alpha\cap[N],\id\circ\phi) \leq \ell(\alpha,\id)+\ell([N],\phi) = \ell_+(\phi)\,,$$
because $\ell(\alpha,\id) = 0$, the identity map being generated by the zero Hamiltonian.

Point (xi) is a consequence of the triangle inequality. For instance,
$$\ell(\alpha,\psi\phi\psi^{-1}) \leq \ell(\alpha,\phi)+\ell_+(\psi)+\ell_+(\psi^{-1}) = \ell(\alpha,\phi)+\ell_+(\psi) - \ell_-(\psi)\,.\qed$$
\end{prf}

\section{Proofs}\label{section_proofs}

\subsection{Main result}

\begin{prf}[of theorem \ref{thm_main_result}]We define $\mu_0 \fc \cG \to \R$ by
$$\mu_0(\phi)=\lim_{k\to \infty}\frac{\ell_+(\phi^k)}{k}\,.$$
The limit exists because the sequence $\{\ell_+(\phi^k)\}_k$ is subadditive. It is finite because of property (ii) in theorem \ref{thm_summary_spectral_invts}. As an immediate consequence of this definition, we obtain point (i).

(ii) Point (viii) in theorem \ref{thm_summary_spectral_invts} implies that for any $\phi,\psi\in\cG$ we have
$$|\ell_+\big((\psi\phi\psi^{-1})^k\big) - \ell_+(\phi^k)|= |\ell_+(\psi\phi^k\psi^{-1}) - \ell_+(\phi^k)| \leq \ell_+(\psi) - \ell_-(\psi)\,.$$
Dividing by $k$ and taking $k \to \infty$ yields $\mu_0(\psi\phi\psi^{-1})=\mu_0(\phi)$.

Let us define the functionals $\mu_a$. For $a \in H^1(N;\R)$ let $\alpha\in a$ and define the symplectomorphism $T_{\alpha} \fc T^*N \to T^*N$ by $T_{\alpha}(q,p)=(q,p+\alpha(q))$. Put
$$\mu_a(\phi)=\mu_0(T_{-\alpha}\phi T_\alpha)\,.$$
Since $\mu_0$ is conjugation-invariant, it follows that if we replace $\alpha$ by a cohomologous form $\alpha'$, then $\mu_0(T_{-\alpha}\phi T_\alpha) = \mu_0(T_{-\alpha'}\phi T_{\alpha'})$ and therefore $\mu_a$ is unambiguously defined. Indeed, assume that $\alpha-\alpha'=df$ for some $f \in C^\infty(N)$. Let $B$ be a cotangent ball large enough to contain the support of a Hamiltonian generating $T_{-\alpha'}\phi T_{\alpha'}$. Let $F$ be a compactly supported Hamiltonian obtained from $\pi^*f$ by cutting it off outside $B$. We then have $T_{-df}T_{-\alpha'}\phi T_{\alpha'}T_{df}=\phi_FT_{-\alpha'}\phi T_{\alpha'}\phi_{-F}$. Consequently
$$T_{-\alpha}\phi T_{\alpha}=T_{-df}T_{-\alpha'}\phi T_{\alpha'}T_{df}=\phi_FT_{-\alpha'}\phi T_{\alpha'}\phi_{-F}$$
and so
$$\mu_0(T_{-\alpha}\phi T_{\alpha})=\mu_0(\phi_F(T_{-\alpha'}\phi T_{\alpha'})\phi_{-F})=\mu_0(T_{-\alpha'}\phi T_{\alpha'})\,.$$

(iii) Let us show the upper bound for $\mu_0$, for example; the rest follows similarly. We have
$$\ell_+(\phi)-\ell_+(\psi) = \ell_+(H)-\ell_+(K) \leq \int_0^1\max(H_t - K_t)\,dt\,.$$
In order to pass to $\mu_0$, we need to homogenize and to this end, to concatenate Hamiltonians. Let $f$ be a smoothing function as in remark \ref{rem_smoothing}. Then
$$\ell_+(H)-\ell_+(K) = \ell_+(H^f)-\ell_+(K^f) \leq \int_0^1\max(H^f_t - K^f_t)\,dt\,.$$
For any $\ve > 0$ there is such a smoothing function for which
$$\int_0^1\max(H^f_t - K^f_t)\,dt \leq \int_0^1\max(H_t - K_t)\,dt + \ve\,.$$
It follows that
$$\frac{\ell_+(\phi^k)-\ell_+(\psi^k)}{k} = \frac{\ell_+\big((H^f)^{\sharp k}\big)-\ell_+\big((K^f)^{\sharp k}\big)}k \leq \int_0^1\max(H_t - K_t)\,dt + \ve\,,$$
and passing to the limit $k\to\infty$, and then letting $\ve \to 0$, we obtain
$$\mu_0(\phi)-\mu_0(\psi) \leq \int_0^1\max(H_t - K_t)\,dt\,.$$

(iv) If $U$ is displaceable by $\psi \in \cG$, we have, by theorem \ref{thm_summary_spectral_invts}:
$$|\ell_+(\phi)| \leq \Gamma(\psi)$$
for any $\phi\in\cG_U$. Then
$$|\mu_0(\phi)| = \lim_{k\to \infty}\frac{|\ell_+(\phi^k)|}{k} \leq \lim_{k\to \infty}\frac{\Gamma(\psi)}{k} = 0\,.$$
For $a\in H^1(N;\R)$ and $\alpha \in a$ we have that $T_{-\alpha}\phi T_{\alpha} \in \cG_{T_{-\alpha}(U)}$ and that $T_{-\alpha}(U)$ is displaceable by $T_{-\alpha}\psi T_\alpha$.

(v) Let us say for brevity that a diffeomorphism $\phi \in \cG$ is dominated by an open subset $U$ if $\phi$ is generated by a Hamiltonian whose support is dominated by $U$. Let $\phi$ be dominated by one of the elements in $\cU$. The triangle inequality and duality for $\ell_\pm$ implies that for any $\alpha,\beta \in \cG$
$$\ell_-(\alpha)+\ell_+(\beta) \leq \ell_+(\alpha\beta)\leq\ell_+(\alpha)+\ell_+(\beta)\,.$$
Put $\phi_j = \psi^j\phi\psi^{-j}$. Then $\phi_j$ is also dominated by one of the elements in $\cU$. We have
$$(\phi\psi)^k=\phi_0\phi_1\dots\phi_{k-1}\psi^k\,,$$
which implies, using induction and the above inequality, that
$$\sum_{j=0}^{k-1}\ell_-(\phi_j)+\ell_+(\psi^k) \leq \ell_+((\phi\psi)^k) \leq \sum_{j=0}^{k-1}\ell_+(\phi_j) + \ell_+(\psi^k)\,.$$
Since all $\phi_j$ are dominated by a displaceable subset with displacement energy $\leq e(\cU)$, it follows, using property (vii) in theorem \ref{thm_summary_spectral_invts}, that
$$|\ell_+((\phi\psi)^k)-\ell_+(\psi^k)|\leq ke(\cU)\,,$$
and upon dividing by $k$ and taking $k\to \infty$ we get
$$|\mu_0(\phi\psi)-\mu_0(\psi)| \leq e(\cU)\,.$$
The claim follows by induction on $\|\phi\|_\cU$.

(vi) Since the restriction of any spectral invariant to $\cG_0$ coincides with the action homomorphism, it follows by homogenization that so does the restriction of $\mu_0$.

(vii) It suffices to restrict the attention to the zero section and $a=0$. From theorem \ref{thm_summary_spectral_invts}, point (iv), we know that if, for example, $H|_N \geq c$, then
$$\ell_+(\phi) \geq c\,.$$
Again, to prove the corresponding property for $\mu_0$, we need to concatenate. Consider the Hamiltonian $K$ given by
$$K_t = c(1-f'(t))+f'(t)H(f(t),\cdot)\,,$$
where $f$ is a smoothing function. Of course, $K$ is not compactly supported, but this is easily circumvented by cutting it off outside a large compact. The action spectrum of $K$ is that of $H^f$, shifted upward by the amount $\int_0^1 c(1-f'(t))\,dt$. This number can be made as small as we wish by suitably choosing $f$. It follows that for any $\ve > 0$ there is a smoothing function $f$ such that the spectral invariants of $K$ differ from those of $H^f$ by not more than $\ve$. We already know that the spectral invariants of $H^f$ and of $H$ coincide. This discussion shows that
$$\ell_+(K) \geq c - \ve\,.$$
But $K$ equals $c$ near $t=0,1$ and so can be concatenated with itself to yield a smooth function. It follows that
$$\frac{\ell_+(\phi^k)}{k} = \frac{\ell_+(K^{\sharp k})}{k} \geq c - \ve\,,$$
and passing first to the limit $k\to \infty$ and then taking $\ve \to 0$ we obtain
$$\mu_0(\phi) \geq c\,.$$

(viii) Follows from the triangle inequality for $\ell_+$ and the fact that $(\phi\psi)^k=\phi^k\psi^k$.

(ix) We have
$$|\mu_a(\phi_H)-\mu_b(\phi_H)| = |\mu_0(T_{-\alpha}\phi_H T_\alpha)-\mu_0(T_{-\beta}\phi_H T_\beta)|\,,$$
where $\alpha \in a, \beta \in b$. The right-hand side is bounded from above by
$$\int_0^1\|H_t\circ T_\alpha - H_t\circ T_\beta\|_{C^0}\,dt\,.$$
For any $1$-form $\chi$ on $N$ we have
$$\max_{T^*N}|H_t - H_t\circ T_\chi| \leq |dH_t|(\chi)\,,$$
where we use the notation
$$|dH_t|(\chi) = \max_{(q,p)\in T^*N}\big|\big\langle d_{(q,p)}H_t|_{T^{\text{vert}}_{(q,p)}T^*N},\chi(q)\big\rangle\big|\,,$$
where we identified $T^{\text{vert}}_{(q,p)}T^*N = T^*_qN$ and $\langle\cdot,\cdot\rangle$ is the pairing between $T^*_qN$ and $T_qN$. It follows that
$$|\mu_a(\phi_H)-\mu_b(\phi_H)| \leq |dH|(a-b)\,,$$
where $|dH| \fc H^1(N;\R)\to\R$ is the semi-norm defined by
$$|dH|(a) = \inf_{\alpha\in a}\int_0^1|dH_t|(\alpha)\,dt\,.$$
This means that $a\mapsto\mu_a(\phi)$ is Lipschitz, the Lipschitz constant being replaced by the semi-norm $|dH|$. \qed
\end{prf}

\begin{rem}
The functions $\mu_a$ have been defined via $\mu_0$, which in turn is the homogenization of the spectral invariant $\ell_+$. An equivalent construction of the $\mu_a$ can be achieved as follows. Fix $\alpha \in a$ and let $L_\alpha$ be the graph of $\alpha$ and $\lambda_\alpha = \lambda - \pi^*\alpha$. Then $L_\alpha$ is an exact Lagrangian submanifold of $(T^*N,\lambda_\alpha)$, on which $\lambda_\alpha$ vanishes, and so we can perform the constructions of section \ref{section_spectral_invariants} in the same fashion, with the zero section being replaced by $L_\alpha$. It is easy to see that this construction leads, through the corresponding version of the spectral invariants $\ell_{+,\alpha}$, to the same functions $\mu_\alpha$.
\end{rem}

Proposition \ref{prop_product_for_qm} follows from the product formula for spectral invariants, theorem \ref{thm_product_formula_sp_invts}, homogenization, and shifting by $T_\alpha$ for appropriate $1$-forms $\alpha$.

\begin{prf}[of theorem \ref{thm_properties_zeta}]Points (ii-iv), (vi), (vii) are immediate consequences of the relevant properties of $\mu_a$. Point (i) is proved by invoking the semi-homogeneity of $\mu_a$ to obtain the desired identity first for natural, then rational, and finally, using continuity, for arbitrary $\lambda \geq 0$. Point (v) is proved as in \cite{EPZ_qm_Poisson_br}, carefully keeping track of the constants. \qed
\end{prf}

\subsection{Equivalence to symplectic homogenization on $\T^n$}\label{section_proof_equivalence_to_homogenization}

Here we prove proposition \ref{thm_our_constr_equals_homogenization} which states that our present construction is equivalent to the symplectic homogenization if $N = \T^n$.

\subsubsection{Overview of the proof}\label{section_overview_equiv_to_homogenization}

Before giving the details, let us present an overview of the construction and an intuitively clear argument why the two constructions are equivalent. We use the notation $\ol{T^*\T^n}$ to indicate that the symplectic form is the negative of the usual one.

One starts with a Hamiltonian $H \in C^\infty_c([0,1]\times T^*\T^n)$ and its flow $\phi^t$. The graph of $\phi^t$ is a Lagrangian submanifold $\Gamma_{\phi^t}\subset T^*\T^n \times \ol{T^*\T^n}$, and it is the image of the diagonal $\Delta\equiv\Delta_{T^*\T^n}\subset T^*\T^n \times \ol{T^*\T^n}$ under the Hamiltonian isotopy $\id \times \phi^t$. There is a symplectic covering $\tau \fc T^*\Delta \to T^*\T^n \times \ol{T^*\T^n}$ which sends the zero section diffeomorphically onto $\Delta$. The isotopy $\id \times \phi^t$ lifts to a unique Hamiltonian isotopy (which no longer has compact support) $\widetilde \Phi^t$ which maps the zero section to a Lagrangian submanifold $L(t) = \widetilde \Phi^t(\cO_\Delta)$. This Lagrangian submanifold maps diffeomorphically onto $\Gamma_{\phi^t}$ under the covering $\tau$.

Since $\phi^t$ has compact support, the Lagrangian $L(t)$ coincides with the zero section outside a compact subset $K$. Consequently it admits a generating function quadratic at infinity $S(t) \fc T^*\Delta \times E \to \R$ ($E$ is a parameter space), which up to a gauge transformation and stabilization is uniquely determined by the requirement that it coincide with a quadratic form on $E$ on the complement of $K\times E$. This implies that its spectral invariants are uniquely determined by $L(t)$, and thus by $\phi(t)$. For $k \in \N$ one defines the function $h_k \fc \R^n \to \R$ by $h_k(p) = \frac 1 k\ell_+(S(k)_p)$ where $S(k)_p = S(k)|_{\T^n \times \{p\}\times E}$, $\T^n \times \{p\}$ being considered a subset of $\Delta = T^*\T^n$.

In fact, Viterbo uses another definition of $h_k$. We will now describe it and show the equality of the two definitions. He uses the covering map $\rho_k \fc T^*\T^n \to T^*\T^n$, $\rho_k(q,p)=(kq,p)$. Being a conformally symplectic covering, this map allows to pull-back Hamiltonian vector fields via the formula $\rho_k^*(X)(z) = (d_z\rho_k)^{-1}(X(\rho_k(z)))$, and thus defines a homomorphism $\rho_k^* \fc \cG \to \cG$. This map enters his construction as follows. Let $H_k \in C^\infty_c([0,1]\times T^*N)$ be defined as\footnote{Formally speaking, one should assume $H$ to be time-periodic for this to make sense, but we suppress such considerations below.} $H_k(t,q,p)=H(kt,kq,p)$. The the time-$1$ flow of $H_k$ is given by $\phi_k:=\phi_{H_k}=\rho_k^*\phi_H^k$.

Now, Viterbo constructs a generating function $S_k$ for the image of the zero section in $T^*\Delta$ under the lift of $\id \times \phi_k$ to $T^*\Delta$ via the above symplectic covering $\tau$. He then defines the function $h_k \fc \R^n \to\R$ by $h_k(p) = \ell_+(S_k|_{\T^n\times\{p\}\times E})$. The two definitions of the function $h_k$ coincide. Indeed, the proof of the next lemma is an easy exercise (see \cite{Viterbo_homogenization}).
\begin{lemma}Let $\phi\in\cG$ and $\psi = \rho_k^*\phi$. Let $S \fc \Delta \times E \to \R$ be the gfqi for the lift of the graph $\Gamma_\phi$ to $T^*\Delta$ as above. Then $T\fc \Delta \times E \to \R$, defined by $T(q,p,\xi)=\frac 1 k S(kq,p,\xi)$ is a gfqi for the lift to $T^*\Delta$ of the graph $\Gamma_\psi$. \qed
\end{lemma}
\noindent It follows that the spectral invariants of $T$ are $\frac 1 k$ times those of $S$. Let us say for brevity that a gfqi for the lift to $T^*\Delta$ via $\tau$ of the graph $\Gamma_\phi$ is a gfqi for $\phi$. Then $S(k)$ is a gfqi for $\phi_H^k$, while $S_k$ is a gfqi for $\phi_k = \rho_k^*\phi_H^k$. It follows from these considerations that the spectral invariants of $S_k$ are $\frac 1 k$ times those of $S(k)$, in particular,
$$\ell_+(S_k|_{\T^n \times\{p\}\times E}) = \frac 1 k \ell_+(S(k)|_{\T^n \times\{p\}\times E})\,,$$
and the two sides of this equality are precisely the two definitions of the function $h_k(p)$ above.

With this function at hand, one can show, using the fact that the sequence $h_k(p)$ is subadditive\footnote{The subadditivity follows \emph{a posteriori} after we have shown that $h_k(0)=\ell_+(\phi_H^k)/k$.} for a fixed $p$, that the limit $h = \lim_kh_k$ exists, and is in fact a continuous function. We denote $\ol H(q,p)=h(p)$, and this is the symplectic homogenization of $H$. It is clear from the above description that if we define $H_p(t,q,\cdot)=H(t,q,\cdot+p)$, then $\ol{H_p}(0)=\ol H (p)$. Since the analogous property holds for the functions $\mu_a$, it suffices to prove
\begin{prop}\label{prop_hbar_equals_mu_at_zero}$\ol H(0)=\mu_0(\phi_H)$.
\end{prop}
\noindent It suffices to see that, for each $k$,
$$h_k(0)=\frac{\ell_+(\phi_H^k)}k\,,$$
since the two sides of the equality asserted in the proposition are obtained as the respective limits when $k \to \infty$. The proposition thus will follow if we show that
$$\ell_+(S(k)|_{\T^n \times\{0\}\times E}) = \ell_+(\phi_H^k)\,.$$
In this form it is almost obvious. The point is that the restriction $S(k)|_{\T^n \times\{0\}\times E}$ generates a Lagrangian submanifold in $T^*\T^n$, which is the image of $\phi_H^k(\cO_{\T^n})$ under the involution $(q,p)\mapsto (q,-p)$. Since the action functional corresponding to $H^{\sharp k}$ generates the same Lagrangian, its spectral invariants will coincide with those of $S(k)_0$ if the two induce the same function on the Lagrangian. The bulk of the proof below is devoted to showing this fact. What allows to conclude is, roughly speaking, the fact that the action functional corresponding to the lifted isotopy $\widetilde \Phi^k$, as well as the gfqi $S(k)$, generate the same Lagrangian $L(k)$, which is the lift to $T^*\Delta$ via $\tau$ of the graph $\Gamma_{\phi^k}$. But these two are normalized to equal zero at points of $\Delta = T^*\T^n$ outside a large compact. It then follows that their spectral invariants, and in particular those of the reduced functionals $S(k)_0$ and $\cA_{H^{\sharp k}}$, coincide. Note the multiple instances of the use of lemma \ref{lemma_sp_invts_FH_gfqi_coincide}.

\subsubsection{Details}\label{section_details_equiv_to_homogenization}

Let us describe the construction of the symplectic homogenization. Fix $H \in C^\infty_c([0,1]\times T^*\T^n)$ and let $\phi^t \equiv \phi_H^t$ be the isotopy generated by $H_t$, and $\phi=\phi^1$. Also denote $\Phi^t = \id \times \phi^t \fc T^*\T^n \times \ol{T^*\T^n} \to T^*\T^n \times \ol{T^*\T^n}$. The isotopy $\Phi^t$ is generated by the Hamiltonian $K_t = 0 \oplus (-H_t)$, that is $K_t(z,z')=-H_t(z')$ for $z,z' \in T^*\T^n$. Consider the symplectic covering $\tau \fc T^*\Delta \to \T^*\T^n \times \ol{T^*\T^n}$, where $\Delta \equiv \Delta_{T^*\T^n} = T^*\T^n$ is the diagonal, given by $\tau(u,v;U,V) = (u-V,v;u,v-U)$. We have the following commutative diagram:
$$\begin{CD}
T^*\R^n \times \ol{T^*\R^n}@>>> T^*\Delta_{T^*\R^n}\\
@VVV @VVV\\
T^*\T^n \times \ol{T^*\T^n} @<\tau<< T^*\Delta_{T^*\T^n}
\end{CD}\,.$$
Here we view explicitly $\T^n=\R^n/\Z^n$. The left and the right arrows are induced from the quotient maps\footnote{Usually if there is a smooth map $f \fc X \to Y$, there is no natural way of associating a smooth map between the corresponding cotangent bundles, however if this map is a local diffeomorphism, then we get the induced map $f_* \fc T^*X \to T^*Y$ given by $f_*(\alpha) = \alpha \circ (d_xf)^{-1}$ for $\alpha \in T^*_xX$, and it is symplectic: $(f_*)^* \omega^{T^*Y} = \omega^{T^*X}$.} $\R^n \times \R^n \to \T^n \times \T^n$ and $T^*\R^n \to T^*\T^n$. The top map is given by $(q,p;Q,P) \mapsto (Q,p;p-P,Q-q)$.

Consider the Hamiltonian $\widetilde H_t = K_t \circ \tau$ on $T^*\Delta$. It generates a lift $\widetilde \Phi^t$ of $\Phi^t$: $\Phi^t \circ \tau = \tau \circ \widetilde \Phi^t$. Denote $L = \widetilde \Phi^1(\cO_\Delta)$. Then $\tau(L) = \Gamma_{\phi}$.

We can now extract spectral invariants from $\widetilde H_t$. This Hamiltonian is not compactly supported, but for finite $t$ it suffices to cut it off outside a large enough ball, and consider the action functional corresponding to that function. By abuse of notation we denote this modified action functional also by $\cA_{\widetilde H}$. It has the same values on all Hamiltonian arcs starting at the zero section and following the flow $\widetilde \Phi^t$ as the original functional before the cutoff.

Again, the fact that $H_t$ has compact support implies that $L$ differs from the zero section only inside a compact subset of $T^*\Delta$ and so we can compactify all the objects in sight to $T^*(\T^n \times S^n) = T^*(\Delta \cup \T^n \times \{\infty\})$. We denote them by the same letters as their counterparts on $T^*\Delta$.

Viterbo gives a formula for a gfqi generating the Lagrangian $L$ \cite{Viterbo_homogenization}. The precise formula is irrelevant, since, as we mentioned in subsection \ref{section_overview_equiv_to_homogenization}, the spectral invariants of a gfqi for $L$ are uniquely determined as soon as we normalize it to equal a quadratic form outside $K\times E$, $K \subset \Delta$ being a certain compact subset. We denote the gfqi of $L$, normalized in this fashion, by $S \fc \Delta\times E\to\R$, until the end of this subsection.

From the definition of $\widetilde H$ it is clear that the points of $\T^n \times \{\infty\}$, considered as constant curves, are Hamiltonian arcs with respect to $\widetilde H_t$, starting and ending at the zero section, and moreover that $\widetilde H_t$ actually equals zero on an open neighborhood of $\T^n \times \{\infty\}$ inside $T^*(\T^n \times S^n)$, which means, in particular, that the action of a point in $\T^n \times \{\infty\}$, considered as a Hamiltonian arc, is zero.

Recall that both $\ol\cA_{\widetilde H}$ and $S$ generate $L$, and so induce functions on $L \subset T^*(\T^n \times S^n)$ (see subsection \ref{section_Lagr_sp_invts_gfs}) and that these functions differ by a constant. Now, it follows from the previous paragraphs that the values of both these functions at a point of $\T^n \times \{\infty\}$ is zero, which implies that they functions coincide. In particular, if $\gamma \fc [0,1] \to T^*(\T^n \times S^n)$ is a Hamiltonian arc relative to $\widetilde H_t$, starting at the zero section, and $z = \gamma(1) \in L$, then
$$(S|L)(z) = \ol\cA_{\widetilde H}(\gamma) = -\cA_{\widetilde H}(\gamma)\,.$$

Symplectic homogenization is defined in terms of the spectral invariants of the functions $S_p \fc \T^n \times E \to \R$, $p \in \R^n$, where $S \fc T^*\T^n \times E \to \R$ is the generating function of $L = \widetilde \Phi^1(\cO_\Delta)$ described above, and $S_p(q;\xi):=S(q,p;\xi)$. It turns out that $S_p$ generates the Lagrangian submanifold of $T^*\T^n$ given by
$$L_p=\{(Q(q,p),p-P(q,p))\,|\,q \in \T^n,(Q(q,p),P(q,p))=\phi(q,p)\}\,.$$
This is a simple computation which can be checked using the above commutative diagram. Since we want to prove proposition \ref{prop_hbar_equals_mu_at_zero}, we restrict ourselves to the case $p=0$, so that $S_0$ generates the following Lagrangian:
$$L_0=\{(Q(q,0),-P(q,0))\,|\,q \in \T^n,(Q(q,0),P(q,0))=\phi(q,0)\} = \ol{\phi(\cO_{\T^n})}\,,$$
where for a Lagrangian $Y \subset T^*X$ we denote by $\ol Y$ the flipped Lagrangian, that is, the image of $Y$ by the involution $(q,p)\mapsto (q,-p)$.

The same Lagrangian submanifold is generated by the action functional $\cA_H$. We want to show that the spectral invariants of $H$ coincide with those of $S_0$, namely $\ell_\pm (H)=\ell_\pm(S_0)$. First, we have
\begin{lemma}\label{lemma_action_gfqi_coincide}$\cA_H|L_0 = S_0|L_0$.
\end{lemma}
\noindent Postponing the proof of the lemma for now, let us show how it allows to conclude. Since both $\cA_H$ and $S$ generate $L_0$, it follows that both $\ol\cA_H = -\cA_H$ and $-S_0$ generate $\ol{L_0} = \phi(\cO_{\T^n})$, and that $(-\cA_H)|\ol{L_0} = (-S_0)|\ol{L_0}$, which yields (see subsections \ref{section_sign_conventions}, \ref{section_relation_FH_gf_invts})
$$\ell_\pm(\phi)=\ell_\pm(H)=-\ol\ell_\mp(H)=-\ell_\mp(-S_0)\,.$$
By duality considerations (see \cite{Viterbo_gfqi} for example),
$$\ell_\pm(-S_0) = -\ell_\mp(S_0)\,,$$
so we obtain finally
$$\ell_\pm(\phi)=\ell_\pm(S_0)\,.$$
Note that the whole construction up to this point can be performed with $\phi$ replaced by $\phi^k$ and therefore
$$\ell_+(\phi^k) = \ell_+(S(k)_0)\,,$$
where, as before, $S(k)$ is the gfqi of the Lagrangian $\widetilde \Phi^k(\cO_\Delta) \subset T^*\Delta$. The discussion after the formulation of proposition \ref{prop_hbar_equals_mu_at_zero} shows that this suffices to prove the proposition and therefore we are done.

It only remains to prove lemma \ref{lemma_action_gfqi_coincide}.
\begin{prf}
Since both $\cA_H$ and $S_0$ generate the same Lagrangian $L_0$, it suffices to show their equality at one point of $L_0$. Choose a point $z \in L_0 \cap \cO_{\T^n}$. It exists by Lagrangian intersection theory. Let $\gamma$ be the Hamiltonian arc ending at $z$, relative to the flow $\phi_H^t$, that is $\gamma(t) = \phi_H^t(\gamma(0))$ and $\gamma(1) = z$. Denote $(q,0) = \gamma(0) \in T^*\T^n$. In coordinates, $\gamma(t)=(Q_t,P_t)$. Note that the curve $t \mapsto Q_t \in \T^n$ has lifts to $\R^n$, and for any such lift, say, $\delta(t)$, the difference $\delta(t)-\delta(0)$ is independent of the lift. We denote this difference by $Q_t-q \in \R^n$.

Consider the following arc $\widetilde \gamma \fc [0,1] \to T^*\Delta$:
$$\widetilde \gamma (t) = \widetilde\Phi^t(\gamma(0))$$
where $\gamma(0) \in T^*\T^n = \Delta$; we have
$$(\tau\circ\widetilde\gamma)(t)=(q,0;\gamma(t))=(q,0;Q_t,P_t)\in T^*\T^n\times\ol{T^*\T^n}$$
and
$$\widetilde \gamma(t) = (Q_t,0;-P_t,Q_t-q) \in T^*\Delta\,.$$
Let us compute the action of this arc relative to the Hamiltonian $\widetilde H_t$, that is,
$$\cA_{\widetilde H}(\widetilde \gamma) = \int_0^1 \widetilde H_t(\widetilde \gamma(t))\,dt - \int \widetilde \gamma^*\lambda^\Delta\,,$$
where $\lambda^\Delta$ is the Liouville form on $T^*\Delta$. We have in the first integral:
$$\int_0^1 \widetilde H_t(\widetilde \gamma(t))\,dt = \int_0^1 (K_t\circ\tau\circ\widetilde\gamma)(t)= \int_0^1 (0\oplus -H_t)(q,0;\gamma(t))\,dt = -\int_0^1 H_t(\gamma(t))\,dt\,.$$
The second integral equals
$$-\int_0^1 \lambda^\Delta(\dot{\widetilde \gamma}(t))\,dt = -\int_0^1 \langle (-P_t,Q_t-q),{\textstyle\frac{d}{dt}(Q_t-q,0)}\rangle\,dt = \int_0^1\langle P_t,\dot Q_t\rangle\,dt = \int \gamma^*\lambda\,.$$
In total we get
$$-\cA_{\widetilde H}(\widetilde \gamma) = \cA_H(\gamma)\,.$$
Denoting $\widetilde z = \widetilde\gamma(1)$, we have
$$(S_0|L_0)(z) = (S|L)(\widetilde z) = -\cA_{\widetilde H}(\widetilde\gamma) = \cA_H(\gamma)\,,$$
as asserted. The first of these equalities follows from the fact that $L_0$ is obtained from $L$ by symplectic reduction (which is just a reformulation of the fact that $L_0$ is generated by the gfqi $S_0$ which itself is the restriction of $S$ to the zero section $\cO_{\T^n} \subset T^*\T^n = \Delta$), and that under this reduction $\widetilde z$ is mapped to $z$. \qed
\end{prf}

\subsection{Alpha function}

\begin{prf}[of theorem \ref{thm_homogen_equals_alpha}]It suffices to show the equality for $a=0$. We have the following expression for the alpha function at zero (this is implicit in \cite{Mather_action_minimizing}; see for example the proof of the proposition on page 178):
$$\alpha_H(0) = -\lim_{k\to\infty}\frac 1 k \inf\{\cA_L^k(\gamma)\,|\,\gamma\fc[0,k]\to N\}\,,$$
where $k \in \N$ and
$$\cA_L^k(\gamma) = \int_0^k L(t,\gamma(t),\dot\gamma(t))\,dt\,,$$
$L \fc \R \times TN \to \R$ being the time-periodic Lagrangian function associated to $H$ by the Fenchel duality. We claim that the infimum in the right-hand side equals $-\ell_+(\phi_H^k)$. Assuming this claim for a moment, we obtain
$$\alpha_H(0) = \lim_{k \to \infty}\frac{\ell_+(\phi_H^k)}{k} = \mu_0(\phi_H)\,,$$
as required.

To prove the claim, consider the functional $\cA_L^k \fc \cP_k \to \R$, where $\cP_k = \{\gamma\fc[0,k]\to N\}$. The evaluation map $\pi_k \fc \cP_k \to N$, $\gamma \mapsto \gamma(k)$, is a submersion, therefore one can consider $\cA_L^k$ as a generating function. It generates the Lagrangian submanifold $\phi_H^k(N)$. The above infimum is in fact a minimum, therefore a critical value of $\cA_L^k$, and as such, it is the action of a Hamiltonian arc running from the zero section back to itself. We would like to show that this critical value is a spectral invariant of $H$. First, it is possible to find a genuine finite-dimensional generating function $S_k$ for $\phi_H^k(N)$ whose associated quadratic form is positive-definite;\footnote{A proof of this fact can be found in the latest version of \cite{Viterbo_homogenization}, appendix D.} in this case
$$\min S_k = \ell_-(S_k)\,.$$
Since any two generating functions for the same Lagrangian submanifold induce functions on it whose difference is constant, by normalizing $S_k$ we can assume that its critical values coincide with those of $\cA_L^k$. Thus
$$\min \cA_L^k = \min S_k = \ell_-(S_k)\,.$$
Our sign conventions imply that the Hamiltonian action functional is the negative of the Lagrangian one when evaluated at a critical point. Therefore $S_k|\phi_H^k(N) = \cA_L^k|\phi_H^k(N) = -\cA_H^k|\phi_H^k(N)$, and it follows that the spectral invariants of $S_k$ coincide with those of $-\cA_H^k$. Therefore, by duality,
$$\min \cA_L^k=\ell_-(S_k) = -\ell_+(\phi_H^k)\,,$$
as claimed. \qed
\end{prf}

\subsection{Hofer geometry and the spectral metric}

\begin{prf}[of theorem \ref{thm_Hofer_embeddings}]Let us first deal with (ii). Let $H \in C^\infty_c(T^*N)$ be such that $H|_N =1$ and $0\leq H\leq 1$ everywhere. Define $\iota \fc \R \to \cG$ by $t \mapsto \phi_{tH}$. We have
$$\rho(\iota(t),\iota(t')) \leq \osc (t-t')H = |t-t'|\,.$$
On the other hand
$$\rho(\iota(t),\iota(t')) \geq |\mu_0(\phi_{tH}) - \mu_0(\phi_{t'H})|=|t-t'|\,,$$
since $tH|_N = t$ and $t'H|_N=t'$.

The argument for (i) is an elaboration of this trick. Fix a non-singular closed $1$-form $\alpha$ and let $a=[\alpha]\in H^1(N;\R)$. Let $H \fc T^*N \to \R$ be smooth, such that the restriction of $H$ to the graph of $t\alpha$ equals $t$ for $t\in[0,1]$. This is possible because $\alpha$ has no zeros. Now define a map $C^\infty_c(0,1) \to C^\infty_c(T^*N)$ by $f\mapsto H_f:=f\circ H$. This is a linear map. Define $\iota\fc C^\infty_c(0,1) \to \cG$ by $\iota(f)\equiv\phi_f:=\phi_{H_f}$. This map is a group homomorphism. We have
$$\max H_f = \max f$$
and same for $\min$ and $\osc$. Consequently
$$\rho(\iota(f),\iota(g)) \leq \osc(H_f-H_g) = \osc(f-g)\,.$$
On the other hand, if $F,G$ are time-dependent Hamiltonians generating $\phi_f,\phi_g$, respectively, then
$$\mu_{ta}(\phi_f) - \mu_{ta}(\phi_g) \leq \int_0^1\max(F_t-G_t)\,dt\,,$$
but
$$\mu_{ta}(\phi_f) - \mu_{ta}(\phi_g) = f(t) - g(t)=(f-g)(t)\,,$$
by construction, since $f\circ H$ equals $f(t)$ on the graph of $t\alpha$, and similarly for $g$. It follows that
$$\max(f-g) \leq \int_0^1\max(F_t-G_t)\,dt\,,$$
and similarly
$$\min(f-g) \geq \int_0^1\min(F_t-G_t)\,dt\,.$$
These two inequalities imply that
$$\osc(f-g) \leq \rho(\phi_f,\phi_g)=\rho(\iota(f),\iota(g))\,.$$
The proof for the spectral metric is analogous. \qed
\end{prf}

Next we prove proposition \ref{prop_asym_Hofer_spectral}.
\begin{prf}If $H$ is a Hamiltonian generating $\phi$, then for any $a\in H^1(N;\R)$ we have
$$\int_0^1 \min H_t\,dt\leq\mu_a(\phi) \leq \int_0^1\max H_t\,dt\,.$$
It follows that
$$\int_0^1 \min H_t\,dt\leq \min_{a\in H^1(N;\R)}\mu_a(\phi) \leq \max_{a\in H^1(N;\R)}\mu_a(\phi) \leq \int_0^1\max H_t\,dt\,,$$
which implies
$$\osc_{a\in H^1(N;\R)}\mu_a(\phi) \leq \int_0^1 \osc H_t\,dt\,,$$
and the assertion about the (asymptotic) Hofer metric follows.
For the spectral metric we have the comparison inequality
$$c_-(\phi)\leq\ell_+(\phi)\leq c_+(\phi)\,.$$
The triangle inequality for $c_\pm$ implies that the sequence $\{c_+(\phi^k)\}_{k \geq 1}$ is subadditive, the sequence $\{c_-(\phi^k)\}_{k \geq 1}$ is superadditive, which means
$$c_-(\phi)\leq\frac{\ell_+(\phi^k)}{k}\leq c_+(\phi)\,,$$
therefore
$$c_-(\phi) \leq \mu_0(\phi) \leq c_+(\phi)\,.$$
The spectral invariants $c_\pm$ are invariant under the symplectomorphisms $T_\alpha$ (see the proof of theorem \ref{thm_main_result} for their definition), therefore
$$\osc_{a\in H^1(N;\R)}\mu_a(\phi) \leq c_+(\phi)-c_-(\phi) = \Gamma(\phi)\,,$$
as claimed. Finally, note that the spectral norm satisfies $\Gamma(\phi) \leq \rho(\phi)$. \qed
\end{prf}

\begin{prf}[of theorem \ref{thm_alpha_fcn_Hofer_metric}]Consider a smooth function $f \fc [0,\infty)\to[0,1]$ such that $f(t) = t$ for $t\in[0,1/2]$, $f(t) = 1$ for $t \geq 2$ and such that $f'(t) \geq 0$ everywhere. For $\ve > 0$ put $f_\ve(t)=\ve f(t/\ve)$. Now define $K_\ve = f_\ve \circ \widetilde H$. Note that the flow of $K_\ve$ is generated by the compactly supported Hamiltonian $K_\ve - \ve$.

We have (see, for example \cite{Sorrentino_Viterbo_act_min_pties_dist_Ham}):
$$\lim_{\ve \to 0}\rho(\phi_{K_\ve}) = \rho_\cH(\phi_H)\,;$$
next, for any $a\in H^1(N;\R)$ such that $\|a\|<1$, it is true that
$$\alpha_{K_\ve}(a)=\alpha_{\widetilde H}(a)\,.$$
Finally, note that the minima $\min \alpha_{\widetilde H}$, $\min \alpha_{K_\ve}$ are both negative and attained on $\{\|a\|<1\} \subset H^1(N;\R)$.\footnote{Here $\|\cdot\|$ is the Gromov-Federer stable norm; see \cite{Paternain_P_S_boundary_rigidity} for more information.} For any $\|a\|<1$ we have
$$\rho(\phi_{K_\ve}) \geq -\mu_a(\phi_{K_\ve}) = -\alpha_{K_\ve}(a)$$
therefore
$$\rho(\phi_{K_\ve}) \geq -\min\alpha_{K_\ve} = -\min\alpha_{\widetilde H}\,,$$
and taking $\ve \to 0$, we obtain the desired inequality. \qed
\end{prf}

\subsection{Symplectic rigidity}\label{section_proofs_symp_rigidity}

\begin{prf}[of lemma \ref{lemma_visible_iff_superheavy}]Assume that $X$ is superheavy and for smooth $f$ let $c = \max_X f$. It is possible, for any $\ve>0$, to find $g \in C^\infty_c(T^*N)$ such that $g|_X=c$ and $\|g-f\|_{C^0} \leq \ve$. Then $\zeta(f)-\zeta(g) \leq \max(f-g)\leq \ve$; on the other hand $\zeta(g) = c$, thus $\zeta(f) \leq c + \ve$. Now take $\ve \to 0$.

Conversely, if $f|_X=c$, then $\zeta(f) \leq c$. On the other hand
$$c =\min_X f = -\max_X(-f) \leq - \zeta(-f) \leq \zeta(f)\,,$$
where we used propertiy (vii) of theorem \ref{thm_properties_zeta}. Therefore $\zeta(f) = c$. \qed
\end{prf}

\begin{prf}[of proposition \ref{prop_examples_of_visible_sets}]We need to show that if $f \in C^\infty_c(T^*N)$ satisfies $f|_X=c$ then $\zeta(f) = c$. By the $C^0$ continuity of $\zeta$ it suffices to show this for all $f$ which equal $c$ on an open neighborhood of $X$. Therefore let $f$ be such a function. Let $\widehat X = X \cup \bigcup_i U_i = T^*N - U_\infty$ and let $\widehat f$ be defined as follows: it coincides with $f$ on $U_\infty$ and equals $c$ on $\widehat X$. Then $\widehat f$ is smooth and $\widehat f|_N = c$ since $\widehat X \supset N$. It follows that $\zeta(\widehat f) = c$. On the other hand, if we define the function $f_i$ by $f_i|_{U_i^c}=0$ and $f_i|_{U_i}=c-f$ for each $i$, it follows that $f_i$ is a smooth function with compact support inside $U_i$, which is displaceable, that all the $f_i$ commute with each other and with $f$, and that $\widehat f = f + \sum_i f_i$. This implies, together with the properties of $\zeta$, that $\zeta (f) = \zeta(\widehat f) = c$. \qed
\end{prf}

\begin{prf}[of theorem \ref{thm_product_visible}]Put $X = X_1 \times X_2$. It is enough to show that for any $f$ such that $f|_X=c$, we have $\zeta(f) \leq c$. Due to the Lipschitz continuity of $\zeta$, it suffices to show the above for any function $f$ which equals $c$ on a neighborhood of $X$. So choose such a function $f$ and let $U$ be the neighborhood. Let $U_i \supset X_i$ be neighborhoods of $X_i$ such that $\ol{U_1\times U_2} \subset U$. Let $S_i$ be a closed cotangent disk bundle in $T^*N_i$ which contains the image under the projection $T^*(N_1\times N_2) \to T^*N_i$ of the support of $f$. Finally, let $M > 0$ be a real number which satisfies $\min(2M, M+c/2) \geq \max f$. Consider functions $f_i\in C^\infty_c(T^*N_i)$ such that $f_i|_{X_i}=c/2$, $f_i|_{U_i - X_i} \geq c/2$, $f_i|_{S_i-U_i^c} = M$, and $f_i|_{S_i^c} \geq 0$. As a case-by-case verification shows, $f_1\oplus f_2 \geq f$ on $S$, and moreover $f_1\oplus f_2$ is positive on a small neighborhood $V$ of $S_1 \times S_2$. Moreover, the flow of $f_1 \oplus f_2$ keeps the zero section inside $S_1 \times S_2$. Let $g$ be a cutoff of $f_1 \oplus f_2$ outside $V$. Then $g$ is a compactly supported function verifying $g \geq f$, and in addition the Hamiltonian flow of $g$ keeps the zero section inside $S_1 \times S_2$. 

The definition of $\zeta$ for Hamiltonians with complete flow shows that if $h$ has complete flow and $h'$ is a cutoff outside a compact which contains the image of the zero section under the flow of $h$, then $\zeta(h) = \zeta(h')$ (lemma \ref{lemma_sp_invts_geom_bdd}). The properties of $\zeta$ (theorem \ref{thm_properties_zeta}) then show that
$$\zeta(f) \leq \zeta(g) = \zeta_1(f_1) + \zeta_2(f_2) = c\,.\qed$$
\end{prf}


\end{document}